# BANDWIDTH SELECTION FOR SMOOTH BACKFITTING IN ADDITIVE MODELS

BY ENNO MAMMEN[1] AND BYEONG U. PARK[2]

*University of Mannheim and Seoul National University*

The smooth backfitting introduced by Mammen, Linton and Nielsen [*Ann. Statist.* **27** (1999) 1443–1490] is a promising technique to fit additive regression models and is known to achieve the oracle efficiency bound. In this paper, we propose and discuss three fully automated bandwidth selection methods for smooth backfitting in additive models. The first one is a penalized least squares approach which is based on higher-order stochastic expansions for the residual sums of squares of the smooth backfitting estimates. The other two are plug-in bandwidth selectors which rely on approximations of the average squared errors and whose utility is restricted to local linear fitting. The large sample properties of these bandwidth selection methods are given. Their finite sample properties are also compared through simulation experiments.

**1. Introduction.** Nonparametric additive models are a powerful technique for high-dimensional data. They avoid the curse of dimensionality and allow for accurate nonparametric estimates also in high-dimensional settings; see Stone [20] among others. On the other hand, the models are very flexible and allow for informative insights on the influences of different covariates on a response variable. This is the reason for the popularity of this approach. Estimation in this model is much more complex than in classical nonparametric regression. Proposed estimates require application of iterative algorithms and the estimates are not given as local weighted sums of independent observations as in classical nonparametric regression. This

Received February 2003; revised July 2004.
[1]Supported by DFG-project MA 1026/6-2 of the Deutsche Forschungsgemeinschaft.
[2]Supported by KOSEF, the Korean–German Cooperative Science Program 2001. The paper was written during a visit of B. U. Park in Heidelberg, financed by KOSEF and DFG, project 446KOR113/155/0-1.

complicates the asymptotic analysis of the estimate. In this paper we discuss practical implementations for the smooth backfitting algorithm. Smooth backfitting was introduced in [9]. In particular, we will discuss data-adaptive bandwidth selectors for this estimate. We will present asymptotic results for the bandwidth selectors. Our main technical tools are uniform expansions of the smooth backfitting estimate of order $o_P(n^{-1/2})$ that allow us to carry over results from classical nonparametric regression.

There have been three main proposals for fitting additive models: the ordinary backfitting procedure of Buja, Hastie and Tibshirani [1], the marginal integration technique of Linton and Nielsen [8] and the smooth backfitting of Mammen, Linton and Nielson [9]. Some asymptotic statistical properties of the ordinary backfitting have been provided by Opsomer and Ruppert [13] and Opsomer [12]. Ordinary backfitting is not oracle efficient, that is, the estimates of the additive components do not have the same asymptotic properties as if the other components were known. The marginal integration estimate is based on marginal integration of a full dimensional regression estimate. The statistical analysis of marginal integration is much simpler. In [8] it is shown for an additive model with two additive components that marginal integration achieves the one-dimensional $n^{-2/5}$ rate of convergence under the smoothness condition that the component functions have two continuous derivatives. However, marginal integration does not produce rate-optimal estimates unless smoothness of the regression function increases with the number of additive components. The smooth backfitting method does not have these drawbacks. It is rate-optimal and its implementation based on local linear estimation achieves the same bias and variance as the oracle estimator, that is, the theoretical estimate that is based on knowing other components. It employs a projection interpretation of popular kernel estimators provided by Mammen, Marron, Turlach and Wand [10], and it is based on iterative calculations of fits to the additive components. A short description of smooth backfitting will be given in the next two sections. This will be done for Nadaraya–Watson kernel smoothing and for local linear fits.

For one-dimensional response variables $Y^i$ and $d$-dimensional covariates $X^i = (X_1^i, \dots, X_d^i)$ $(i = 1, \dots, n)$ the additive regression model is defined as

$$Y^i = m_0 + \sum_{j=1}^{d} m_j(X_j^i) + \varepsilon^i, \tag{1.1}$$

where $X^i = (X_1^i, \dots, X_d^i)$ are random design points in $\mathbf{R}^d$, $\varepsilon^i$ are unobserved error variables, $m_1, \dots, m_d$ are functions from $\mathbf{R}$ to $\mathbf{R}$ and $m_0$ is a constant. Throughout the paper we will make the assumption the tuples $(X^i, \varepsilon^i)$ are i.i.d. and that the error variables $\varepsilon^i$ have conditional mean zero (given the covariates $X^i$). Furthermore, it is assumed that $Em_j(X_j^i) = 0$ for $j = 1, \dots, d$ and that $\sum_{j=1}^d f_j(X_j^i) = 0$ a.s. implies $f_j \equiv 0$ for all $j$. Then the functions $m_j$ are uniquely identifiable. The latter assumption is a sufficient condition to avoid *concurvity* as termed by Hastie and Tibshirani [5].

Our main results are higher-order stochastic expansions for the residual sums of squares of the smooth backfitting estimates. These results motivate the definition of a penalized sum of squared residuals. The bandwidth that minimizes the penalized sum will be called *penalized least squares bandwidth*. We will compare the penalized sum of squares with the average weighted squared error ($ASE$)

$$ASE = n^{-1} \sum_{i=1}^{n} w(X^i) \left\{ \widehat{m}_0 + \sum_{j=1}^{d} \widehat{m}_j(X_j^i) - m_0 - \sum_{j=1}^{d} m_j(X_j^i) \right\}^2. \tag{1.2}$$

Here $w$ is a weight function. We will show that up to an additive term which is independent of the bandwidth the average weighted squared error is asymptotically equivalent to the penalized sum of squared residuals. This implies that the penalized least squares bandwidth is asymptotically optimal. The results for Nadaraya–Watson smoothing are given in the next section. Local linear smoothing will be discussed in Section 3.

In addition to the penalized least squares bandwidth choice, we discuss two plug-in selectors. The first of these is based on a first-order expansion of $ASE$ given in (1.2). This error criterion measures accuracy of the sum of the additive components. An alternative error criterion measures the accuracy of each single additive component,

$$ASE_j = n^{-1} \sum_{i=1}^{n} w_j(X_j^i) \{\widehat{m}_j(X_j^i) - m_j(X_j^i)\}^2.$$

Here $w_j$ is a weight function. Use of $ASE_j$ instead of $ASE$ may be motivated by a more data-analytic focus of the statistical analysis. Additionally, a more technical advantage holds for local linear smoothing. The first-order expansion of $ASE_j$ only depends on the corresponding single bandwidth and does not involve the bandwidths of the other components. In particular, the plug-in bandwidth selector based on the approximation of $ASE_j$ can be written down explicitly. For Nadaraya–Watson backfitting estimates the bias of a single additive component depends on the whole vector of bandwidths. Therefore an asymptotic expansion of $ASE_j$ involves the bandwidths of all components. Also for the global error criterion $ASE$ implementation of plug-in rules for Nadaraya–Watson smoothing is much more complex. The bias part in the expansion of $ASE$ for the Nadaraya–Watson smoothing has terms related to the multivariate design density, a well-known fact also in the single smoother case, and the bias expression may not even be expressed in a closed form. For these reasons, our discussion on plug-in bandwidths will be restricted to local linear fits.

In classical nonparametric regression, the penalized sum of squared residuals which we introduce in this paper is asymptotically equivalent to cross-validation [4]. We conjecture that the same holds for additive models. The approach based on penalized sum of squared residuals is computationally more feasible than cross-validation. It only requires one $n$th of the computing time that is needed for the latter. In the numerical study presented in Section 5, we found that the penalized least squares bandwidth is a good approximation of the stochastic $ASE$-minimizer. It turned out that it outperforms the two plug-in bandwidths by producing the least $ASE$, while for accuracy of each one-dimensional component estimator, that is, in terms of $ASE_j$, none of the bandwidth selectors dominates the others in all cases. In general, plug-in bandwidth selection requires estimation of additional functionals of the regression function (and of the design density). For this estimation one needs to select other tuning constants or bandwidths. Quantification of the optimal secondary tuning constant needs further asymptotic analysis and it would require more smoothness assumptions on the regression and density functions. See [15], [16] and [19]. In this paper, we do not pursue this issue for the plug-in selectors. We only consider a simple choice of the auxiliary bandwidth.

In this paper we do not address bandwidth choice under model misspecification. For additive models this is an important issue because in many applications the additive model will only be assumed to be a good approximation for the true model. We conjecture that the penalized least squares bandwidth will work reliably also under misspecification of the additive model. This conjecture is supported by the definition of this bandwidth. Performance of the plug-in rules has to be carefully checked because in their definitions they make use of the validity of the additive model.

There have been many proposals for bandwidth selection in density and regression estimation with single smoothers. See [17] and [7] for kernel density estimation, and [6] for kernel regression estimation. For additive models there have been only a few attempts for bandwidth selection. These include [14] where a plug-in bandwidth selector is proposed for the ordinary backfitting procedure, [21] where generalized cross-validation is applied to penalized regression splines and [11] where cross-validation is discussed for smooth backfitting.

In this paper we discuss smooth backfitting for Nadaraya–Watson smoothing and for local linear smoothing. For practical implementations we definitely recommend application of local linear smoothing. Local linear smooth backfitting achieves oracle bounds. The asymptotic bias and variance of the estimate of an additive component do not depend on the number and shape of the other components. They are the same as in a classical regression model with one component. This does not hold for Nadaraya–Watson smoothing.

Nevertheless in this paper we have included the discussion of Nadaraya–Watson smoothing. This has been done mainly for clarity of exposition of ideas and proofs. Smooth backfitting with local linear smoothing requires a much more involved notation. This complicates the mathematical discussions. For this reason we will give detailed proofs only for Nadaraya–Watson smoothing. Ideas of the proofs carry over to local linear smoothing. In Section 2 we start with Nadaraya–Watson smoothing. Smooth backfitting for local linear smoothing is treated in Section 3. Practical implementations of our bandwidth selectors are discussed in Section 4. In Section 5 simulation results are presented for the performance of the discussed bandwidth selectors. Section 6 states the assumptions and contains the proofs of the theoretical results.

**2. Smooth backfitting with Nadaraya–Watson smoothing.** We now define the smooth Nadaraya–Watson backfitting estimates. The estimate of the component function $m_j$ in (1.1) is denoted by $\widehat{m}_j^{NW}$. We suppose that the covariates $X_j$ take values in a bounded interval $I_j$. The backfitting estimates are defined as the minimizers of the following smoothed sum of squares:

$$\sum_{i=1}^{n} \int_I \left\{ Y^i - \widehat{m}_0^{NW} - \sum_{j=1}^{d} \widehat{m}_j^{NW}(u_j) \right\}^2 K_h(u, X^i)\, du. \tag{2.1}$$

The minimization is done under the constraints

$$\int_{I_j} \widehat{m}_j^{NW}(u_j) \widehat{p}_j(u_j)\, du_j = 0, \qquad j = 1, \ldots, d. \tag{2.2}$$

Here, $I = I_1 \times \cdots \times I_d$ and $K_h(u, x^i) = K_{h_1}(u_1, x_1^i) \cdot \cdots \cdot K_{h_d}(u_d, x_d^i)$ is a $d$-dimensional product kernel with factors $K_{h_j}(u_j, v_j)$ that satisfy for all $v_j \in I_j$

$$\int_{I_j} K_{h_j}(u_j, v_j)\, du_j = 1. \tag{2.3}$$

The kernel $K_{h_j}$ may depend also on $j$. This is suppressed in the notation. In (2.2) $\widehat{p}_j$ denotes the kernel density estimate of the density $p_j$ of $X_j^i$,

$$\widehat{p}_j(u_j) = n^{-1} \sum_{i=1}^{n} K_{h_j}(u_j, X_j^i). \tag{2.4}$$

The usual choice for $K_{h_j}(u_j, v_j)$ with (2.3) is given by

$$K_{h_j}(u_j, v_j) = \frac{K[h_j^{-1}(v_j - u_j)]}{\int_{I_j} K[h_j^{-1}(v_j - w_j)]\, dw_j}. \tag{2.5}$$

Note that for $u_j, v_j$ in the interior of $I_j$ we have

$$K_{h_j}(u_j, v_j) = h_j^{-1} K[h_j^{-1}(v_j - u_j)]$$

when $K$ integrates to 1 on its support.

By differentiation one can show that a minimizer of (2.1) satisfies for $j = 1, \ldots, d$ and $u_j \in I_j$

$$\sum_{i=1}^n \int_{I_{-j}} \left\{ Y^i - \widehat{m}_0^{NW} - \sum_{k=1}^d \widehat{m}_k^{NW}(u_k) \right\} K_h(u, X^i)\, du_{-j} = 0,$$

and thus

$$\sum_{i=1}^n \int_{I} \left\{ Y^i - \widehat{m}_0^{NW} - \sum_{k=1}^d \widehat{m}_k^{NW}(u_k) \right\} K_h(u, X^i)\, du = 0,$$

where $I_{-j} = I_1 \times \cdots \times I_{j-1} \times I_{j+1} \times \cdots \times I_d$ and $u_{-j} = (u_1, \ldots, u_{j-1}, u_{j+1}, \ldots, u_d)$. Now, because of (2.3) we can rewrite these equations as

$$\widehat{m}_j^{NW}(u_j) = \widetilde{m}_j^{NW}(u_j) - \sum_{k \neq j} \int_{I_k} \widehat{m}_k^{NW}(u_k) \frac{\widehat{p}_{jk}(u_j, u_k)}{\widehat{p}_j(u_j)}\, du_k - \widehat{m}_0^{NW}, \tag{2.6}$$

$$\widehat{m}_0^{NW} = n^{-1} \sum_{i=1}^n Y^i, \tag{2.7}$$

where $\widehat{p}_{jk}(u_j, u_k) = n^{-1} \sum_{i=1}^n K_{h_j}(u_j, X_j^i) K_{h_k}(u_k, X_k^i)$ is a two-dimensional kernel density estimate of the marginal density $p_{jk}$ of $(X_j^i, X_k^i)$. Furthermore, $\widetilde{m}_j^{NW}(u_j)$ denotes the Nadaraya–Watson estimate

$$\widetilde{m}_j^{NW}(u_j) = \widehat{p}_j(u_j)^{-1} n^{-1} \sum_{i=1}^n K_{h_j}(u_j, X_j^i) Y^i.$$

In case one does not use kernels that satisfy (2.3), equations (2.6) and (2.7) have to be replaced by slightly more complicated equations; see [9] for details.

Suppose now that $I_j = [0, 1]$, and define for a weight function $w$ and a constant $C'_H > 0$

$$\begin{aligned} RSS(h) = n^{-1} \sum_{i=1}^n \mathbf{1}&(C'_H n^{-1/5} \le X_j^i \le 1 - C'_H n^{-1/5} \text{ for } 1 \le j \le d) \\ &\times w(X^i) \{Y^i - \widehat{m}_0^{NW} - \widehat{m}_1^{NW}(X_1^i) - \cdots - \widehat{m}_d^{NW}(X_d^i)\}^2, \end{aligned} \tag{2.8}$$

$$\begin{aligned} ASE(h) = n^{-1} \sum_{i=1}^n \mathbf{1}&(C'_H n^{-1/5} \le X_j^i \le 1 - C'_H n^{-1/5} \text{ for } 1 \le j \le d) \\ &\times w(X^i) \{\widehat{m}_0^{NW} + \widehat{m}_1^{NW}(X_1^i) + \cdots + \widehat{m}_d^{NW}(X_d^i) \\ &\qquad - m_0 - m_1(X_1^i) - \cdots - m_d(X_d^i)\}^2, \end{aligned} \tag{2.9}$$

where $\mathbf{1}(A)$ denotes the indicator which equals 1 if $A$ occurs and 0 otherwise. The indicator function has been included in (2.8) and (2.9) to exclude boundary regions of the design where the Nadaraya–Watson smoother has bias terms of order $n^{-1/5}$. In the following Theorems 2.1 and 2.2 we will consider bandwidths $h_j$ that are smaller than $C'_H n^{-1/5}$. Because we assume that the kernel $K$ has support $[-1,1]$ [see assumption (A1) in Section 6.1], boundary regions with higher-order bias terms are then excluded. We now state our first main result. The assumptions can be found in Section 6.

THEOREM 2.1. *Suppose that assumptions* (A1)–(A4) *apply for model* (1.1) *and that* $\widehat{m}_j^{NW}$ *are defined according to* (2.1) *and* (2.2). *Assume that* $I_j$ *are bounded intervals* ($I_j = [0,1]$ *w.l.o.g.) and that* $K_{h_j}(u_j, v_j)$ *are kernels that satisfy* $K_{h_j}(u_j, v_j) = h_j^{-1} K[h_j^{-1}(v_j - u_j)]$ *for* $h_j \le v_j \le 1 - h_j$ *for a function* $K$ *and* $K_{h_j}(u_j, v_j) = 0$ *for* $|v_j - u_j| \ge h_j$. *Then with* $C'_H$ *as in* (2.8) *and* (2.9) *and for all constants* $C_H < C'_H$, *we have uniformly for* $C_H n^{-1/5} \le h_j \le C'_H n^{-1/5}$

$$\begin{aligned} RSS(h) - n^{-1} \sum_{i=1}^n \mathbf{1}(C'_H n^{-1/5} \le X_j^i \le 1 - C'_H n^{-1/5} \\ \text{for } 1 \le j \le d) w(X^i)(\varepsilon^i)^2 \\ + 2n^{-1} \left\{ \sum_{i=1}^n w(X^i)(\varepsilon^i)^2 \right\} \left\{ K(0) \sum_{j=1}^d \frac{1}{n h_j} \right\} - ASE(h) = o_p(n^{-4/5}). \end{aligned} \tag{2.10}$$

*Furthermore, for fixed sequences* $h$ *with* $C_H n^{-1/5} \le h_j \le C'_H n^{-1/5}$, *this difference is of order* $O_p(n^{-9/10})$.

To state the second main result, let $\beta_j(h, u_j)$, $j = 1, \ldots, d$, denote minimizers of $\int \{\beta(h,u) - \beta_1(h,u_1) - \cdots - \beta_d(h,u_d)\}^2 p(u)\, du$, where

$$\beta(h,u) = \sum_{j=1}^d \left\{ m'_j(u_j) \frac{\partial \log p}{\partial u_j}(u) + \frac{1}{2} m''_j(u_j) \right\} h_j^2 \int t^2 K(t)\, dt.$$

The functions $\beta_j(h, u_j)$, $j = 1, \ldots, d$, are uniquely defined only up to an additive constant. However, their sum is uniquely defined. Define

$$PLS(h) = RSS(h) \left\{ 1 + 2 \sum_{j=1}^d \frac{1}{n h_j} K(0) \right\}.$$

THEOREM 2.2. *Under the assumptions of Theorem* 2.1, *we have uniformly for* $C_H n^{-1/5} \le h_j \le C'_H n^{-1/5}$,

$$ASE(h) = \frac{1}{n} \sum_{i=1}^n w(X^i)(\varepsilon^i)^2 \int K^2(t)\, dt \sum_{j=1}^d \frac{1}{n h_j}$$

$$+\int_I \left\{ \sum_{j=1}^{d} \beta_j(h, u_j) \right\}^2 w(u)p(u)\,du + o_p(n^{-4/5}), \tag{2.11a}$$

$$\begin{aligned} PLS(h) - ASE(h) \\ &= n^{-1} \sum_{i=1}^{n} \mathbf{1}(C'_H n^{-1/5} \le X_j^i \le 1 - C'_H n^{-1/5} \\ &\qquad \text{for } 1 \le j \le d) w(X^i)(\varepsilon^i)^2 \\ &\quad + o_p(n^{-4/5}). \end{aligned} \tag{2.11b}$$

Now we define

$$\widehat{h}_{PLS} = \arg\min PLS(h),$$
$$\widehat{h}_{ASE} = \arg\min ASE(h).$$

Here and throughout the paper, the "arg min" runs over $h$ with $C_H n^{-1/5} \le h_j \le C'_H n^{-1/5}$. It would be a more useful result to have some theory for a bandwidth selector that estimates the optimal bandwidth over a range of rates, for example, $h_j \in [An^{-a}, Bn^{-b}]$ for some prespecified positive constants $a, b, A, B$. This would involve uniform expansions of $RSS(h)$ and $ASE(h)$ over the extended range of the bandwidth, which undoubtedly makes the derivations much more complicated. Thus, it is avoided in this paper.

The following corollary is an immediate consequence of Theorem 2.2.

COROLLARY 2.3. *Under the conditions of Theorem* 2.1

$$\widehat{h}_{PLS} - \widehat{h}_{ASE} = o_p(n^{-1/5}).$$

We conjecture that $(\widehat{h}_{PLS} - \widehat{h}_{ASE})/\widehat{h}_{ASE} = O_p(n^{-1/10})$. This is suggested by the fact that for fixed $h$ in Theorem 2.2 the error term $o_p(n^{-4/5})$ can be replaced by $O_p(n^{-9/10})$.

## 3. Smooth backfitting using local linear fits.

The smooth backfitting local linear estimates are defined as minimizers of

$$\begin{aligned} \sum_{i=1}^{n} \int_I \Bigg\{ Y^i - \widehat{m}_0^{LL} - \sum_{j=1}^{d} \widehat{m}_j^{LL}(u_j) \\ - \sum_{j=1}^{d} \widehat{m}_j^{LL,1}(u_j)(X_j^i - u_j) \Bigg\}^2 K_h(u, X^i)\,du. \end{aligned} \tag{3.1}$$

Here $\widehat{m}_j^{LL}$ is an estimate of $m_j$ and $\widehat{m}_j^{LL,1}$ is an estimate of its derivative.

By using slightly more complicated arguments than those in Section 6 one can show that $\widehat{m}_0^{LL}, \dots, \widehat{m}_d^{LL,1}$ satisfy the equations

$$\begin{aligned}\begin{pmatrix}\widehat{m}_j^{LL}(u_j)\\ \widehat{m}_j^{LL,1}(u_j)\end{pmatrix} &= -\begin{pmatrix}\widehat{m}_0^{LL}\\ 0\end{pmatrix} + \begin{pmatrix}\widetilde{m}_j^{LL}(u_j)\\ \widetilde{m}_j^{LL,1}(u_j)\end{pmatrix}\\ &\qquad - \widehat{M}_j(u_j)^{-1}\sum_{l\neq j}\int \widehat{S}_{lj}(u_l,u_j)\begin{pmatrix}\widehat{m}_l^{LL}(u_l)\\ \widehat{m}_l^{LL,1}(u_l)\end{pmatrix}du_l,\end{aligned} \tag{3.2}$$

$$\begin{aligned}\widehat{m}_0^{LL} &= n^{-1}\sum_{i=1}^n Y^i - \sum_{j=1}^d \int \widehat{m}_j^{LL}(u_j)\widehat{p}_j(u_j)\,du_j\\ &\qquad - \sum_{j=1}^d \int \widehat{m}_j^{LL,1}(u_j)\widehat{p}_j^1(u_j)\,du_j.\end{aligned} \tag{3.3}$$

Here and below,

$$\widehat{M}_j(u_j) = n^{-1}\sum_{i=1}^n K_{h_j}(u_j, X_j^i)\begin{pmatrix}1 & X_j^i - u_j\\ X_j^i - u_j & (X_j^i - u_j)^2\end{pmatrix},$$

$$\widehat{S}_{lj}(u_l,u_j) = n^{-1}\sum_{i=1}^n K_{h_l}(u_l, X_l^i)K_{h_j}(u_j, X_j^i)\begin{pmatrix}1 & X_l^i - u_l\\ X_j^i - u_j & (X_l^i - u_l)(X_j^i - u_j)\end{pmatrix},$$

$$\widehat{p}_j^1(u_j) = n^{-1}\sum_{i=1}^n K_{h_j}(u_j, X_j^i)(X_j^i - u_j),$$

and $\widehat{p}_j$ is defined as in the last section. For each $j$, the estimates $\widetilde{m}_j^{LL}$ and $\widetilde{m}_j^{LL,1}$ are the local linear fits obtained by regression of $Y^i$ onto $X_j^i$; that is, these quantities minimize

$$\sum_{i=1}^n \{Y^i - \widetilde{m}_j^{LL}(u_j) - \widetilde{m}_j^{LL,1}(u_j)(X_j^i - u_j)\}^2 K_{h_j}(u_j, X_j^i).$$

A detailed discussion on why (3.1) is equivalent to (3.2) and (3.3) can be found in [9], where a slightly different notation was used. The definition of $\widehat{m}_0^{LL}, \dots, \widehat{m}_d^{LL,1}$ can be made unique by imposing the additional norming conditions

$$\int \widehat{m}_j^{LL}(u_j)\widehat{p}_j(u_j)\,du_j + \int \widehat{m}_j^{LL,1}(u_j)\widehat{p}_j^1(u_j)\,du_j = 0. \tag{3.4}$$

The smooth backfitting estimates can be calculated by iterative application of (3.2). In each application the current versions of $\widehat{m}_l^{LL}, \widehat{m}_l^{LL,1}$ ($l \neq j$) are plugged into the right-hand side of (3.2) and are used to update $\widehat{m}_j^{LL}, \widehat{m}_j^{LL,1}$. The iteration converges with geometric rate (see [9]). The number of iterations may be determined by a standard error criterion. After the last iteration, a norming constant can be subtracted from the last fit of $\widehat{m}_j^{LL}$ so that (3.4) holds. Because of (3.4) this yields $\widehat{m}_0^{LL} = n^{-1}\sum_{i=1}^n Y^i$.

We now define the residual sum of squares $RSS(h)$ and the average squared error. This is done similarly as in (2.8) and (2.9). But now the sums run over the full intervals $I_j$. This differs from Nadaraya–Watson smoothing where the summation excludes boundary values. For Nadaraya–Watson smoothing the boundary values are removed because of bias problems. Let

$$RSS(h) = n^{-1}\sum_{i=1}^n w(X^i)\{Y^i - \widehat{m}_0^{LL} - \widehat{m}_1^{LL}(X_1^i) - \cdots - \widehat{m}_d^{LL}(X_d^i)\}^2, \tag{3.5}$$

$$\begin{aligned} ASE(h) = n^{-1}\sum_{i=1}^n w(X^i)\{&\widehat{m}_0^{LL} + \widehat{m}_1^{LL}(X_1^i) + \cdots + \widehat{m}_d^{LL}(X_d^i) \\ &- m_0 - m_1(X_1^i) - \cdots - m_d(X_d^i)\}^2. \end{aligned} \tag{3.6}$$

As in Section 6 we define the penalized sum of squared residuals

$$PLS(h) = RSS(h)\left\{1 + 2\sum_{j=1}^d \frac{1}{nh_j}K(0)\right\}.$$

The penalized least squares bandwidth $\widehat{h}_{PLS}$ is again given by

$$\widehat{h}_{PLS} = \arg\min PLS(h).$$

Define

$$\beta_j(u_j) = \tfrac{1}{2}m_j''(u_j)\int t^2 K(t)\,dt.$$

Analogous to Theorems 2.1, 2.2 and Corollary 2.3, we now get the following results for local linear smoothing.

THEOREM 3.1. *Suppose that assumptions* (A1)–(A4) *apply, that* $I_j = [0,1]$ *and that* $\widehat{m}^{LL}$ *is defined according to* (3.1) *and* (3.4) *with kernels* $K_{h_j}(u_j, v_j)$. *The kernels are supposed to satisfy the conditions of Theorem* 2.1. *Then, uniformly for* $C_H n^{-1/5} \le h_j \le C_H' n^{-1/5}$,

$$\begin{aligned} &RSS(h) - n^{-1}\sum_{i=1}^n w(X^i)(\varepsilon^i)^2 + 2n^{-1}\left\{\sum_{i=1}^n w(X^i)(\varepsilon^i)^2\right\}\left\{K(0)\sum_{j=1}^d \frac{1}{nh_j}\right\} \\ &\quad - ASE(h) = o_p(n^{-4/5}), \end{aligned} \tag{3.7}$$

$$\begin{aligned} ASE(h) &= \left\{\int_I w(u)p(u)E[(\varepsilon^i)^2|X^i=u]\,du\right\}\int K^2(t)\,dt\sum_{j=1}^d \frac{1}{nh_j} \\ &\quad + \int_I\left\{\sum_{j=1}^d h_j^2\beta_j(u_j)\right\}^2 w(u)p(u)\,du + o_p(n^{-4/5}), \end{aligned} \tag{3.8}$$

$$PLS(h) - ASE(h) = \frac{1}{n}\sum_{i=1}^n w(X^i)(\varepsilon^i)^2 + o_p(n^{-4/5}), \tag{3.9}$$

$$\widehat{h}_{PLS} - \widehat{h}_{ASE} = o_p(n^{-1/5}). \tag{3.10}$$

*For fixed sequences $h$ with $C_H n^{-1/5} \le h_j \le C'_H n^{-1/5}$, the expansions in (3.7)–(3.9) hold up to order $O_p(n^{-9/10})$.*

If the errors of the expansions in (3.7)–(3.9) would be of order $O_p(n^{-9/10})$, uniformly in $h$, this would imply $(\widehat{h}_{PLS} - \widehat{h}_{ASE})/\widehat{h}_{ASE} = O_p(n^{-1/10})$.

Next we consider plug-in bandwidth selectors. As for penalized least squares, plug-in bandwidth selectors may be constructed that approximately minimize $ASE(h)$. Let $AASE(h)$ denote the nonstochastic first-order expansion of $ASE(h)$, given in (3.8). Define

$$\begin{aligned} \widehat{AASE}(h) &= n^{-1}\sum_{i=1}^n w(X^i)(\widehat{\varepsilon}^i)^2\left\{\int K^2(t)\,dt\right\}\sum_{j=1}^d\frac{1}{nh_j} \\ &\quad + \frac{1}{4n}\sum_{i=1}^n w(X^i)\left\{\sum_{j=1}^d h_j^2\widehat{m}_j''(X_j^i)\right\}^2\left\{\int t^2K(t)\,dt\right\}^2. \end{aligned}$$

Here $\widehat{m}_j''$ is an estimate of $m_j''$ and $\widehat{\varepsilon}^i = Y^i - \widetilde{m}(X^i)$ are residuals based on an estimate $\widetilde{m}(x)$ of $m_0 + m_1(x_1) + \cdots + m_d(x_d)$. Choices of $\widehat{m}_j''$ and $\widetilde{m}$ will be discussed below. A plug-in bandwidth $\widehat{h}_{PL} = (\widehat{h}_{PL,1},\ldots,\widehat{h}_{PL,d})$ is defined by

$$\widehat{h}_{PL} = \arg\min \widehat{AASE}(h). \tag{3.11}$$

The plug-in bandwidth $\widehat{h}_{PL}$ will be compared with the theoretically optimal bandwidth $h_{\text{opt}}$,

$$h_{\text{opt}} = \arg\min AASE(h). \tag{3.12}$$

There is an alternative way of plug-in bandwidth selection for another error criterion. It is based on an error criterion that measures accuracy of each one-dimensional additive component separately. Let

$$ASE_j(h) = n^{-1}\sum_{i=1}^n w_j(X_j^i)\{\widehat{m}_j^{LL}(X_j^i) - m_j(X_j^i)\}^2, \tag{3.13}$$

where $w_j$ is a smooth weight function. It may be argued that $ASE_j$ is more appropriate if the focus is more data-analytic interpretation of the data whereas use of $ASE$ may be more appropriate for finding good prediction rules. Our next result shows that in first-order $ASE_j(h)$ only depends on $h_j$. This motivates a simple plug-in bandwidth selection rule. An analogous result does not hold for Nadaraya–Watson smoothing.

THEOREM 3.2. *Under the assumptions of Theorem 3.1, it holds that, uniformly for $h$ with $C_H n^{-1/5} \le h_l \le C'_H n^{-1/5}$ $(1 \le l \le d)$,*

$$\begin{aligned} ASE_j(h) = {} & \left\{ \int_{I_j} w_j(u_j) p_j(u_j) E[(\varepsilon^i)^2 | X_j^i = u_j] \, du_j \right\} \left\{ \int K^2(t) \, dt \right\} \frac{1}{n h_j} \\ & + \frac{1}{4} h_j^4 \int_{I_j} m_j''(u_j)^2 w_j(u_j) p_j(u_j) \, du_j \left\{ \int t^2 K(t) \, dt \right\}^2 \\ & + o_p(n^{-4/5}). \end{aligned} \tag{3.14}$$

The first-order expansion of $ASE_j(h)$ in (3.14) is minimized by

$$\begin{aligned} h^*_{\text{opt},j} = {} & n^{-1/5} \left[ \left\{ \int_{I_j} w_j(u_j) p_j(u_j) E[(\varepsilon^i)^2 | X_j^i = u_j] \, du_j \right\} \left\{ \int K^2(t) \, dt \right\} \right]^{1/5} \\ & \times \left[ \int_{I_j} m_j''(u_j)^2 w_j(u_j) p_j(u_j) \, du_j \left\{ \int t^2 K(t) \, dt \right\}^2 \right]^{-1/5}. \end{aligned}$$

We note that $h_{\text{opt}}$ defined in (3.12) is different from $h^*_{\text{opt}} = (h^*_{\text{opt},1}, \dots, h^*_{\text{opt},d})$. Now this bandwidth can be estimated by

$$\begin{aligned} \widehat{h}^*_{PL,j} = {} & n^{-1/5} \left[ n^{-1} \sum_{i=1}^n w_j(X_j^i) (\widehat{\varepsilon}^i)^2 \left\{ \int K^2(t) \, dt \right\} \right]^{1/5} \\ & \times \left[ n^{-1} \sum_{i=1}^n w_j(X_j^i) \widehat{m}_j''(X_j^i)^2 \left\{ \int t^2 K(t) \, dt \right\}^2 \right]^{-1/5} \end{aligned} \tag{3.15}$$

with an estimate $\widehat{m}_j''$ of $m_j''$ and residuals $\widehat{\varepsilon}^i = Y^i - \widetilde{m}(X^i)$ based on an estimate $\widetilde{m}(x)$ of $m_0 + m_1(x_1) + \cdots + m_d(x_d)$. Contrary to $\widehat{h}_{PL}$, approximation of the bandwidth selector $\widehat{h}^*_{PL}$ does not require a grid search on a high-dimensional bandwidth space or an iterative procedure with a one-dimensional grid search. See the discussion at the end of Section 4.

Now we present a procedure for estimating $m_j''$, which is required to implement $\widehat{h}_{PL}$ and $\widehat{h}^*_{PL}$. A simple estimate of $m_j''$ may be given by smoothed differentiation of $\widehat{m}_j^{LL}$. However, a numerical study for this estimate revealed that it suffers from serious boundary effects. We propose to use an alternative estimate which is based on a local quadratic fit. It is defined by

$$\widehat{m}_j''(u_j) = 2\widehat{\beta}_{j,2}(u_j), \tag{3.16}$$

where $\widehat{\beta}_{j,2}(u_j)$ along with $\widehat{\beta}_{j,0}(u_j)$ and $\widehat{\beta}_{j,1}(u_j)$ minimizes

$$\int_{I_j} \{\widehat{m}_j^{LL}(v_j) - \widehat{\beta}_{j,0}(u_j) - \widehat{\beta}_{j,1}(u_j)(v_j - u_j) - \widehat{\beta}_{j,2}(u_j)(v_j - u_j)^2\}^2 \times L[g_j^{-1}(v_j - u_j)]\, dv_j.$$

The definitions of $\widehat{h}_{PL}$ and $\widehat{h}_{PL}^*$ make use of fitted residuals. But these residuals along with the local quadratic estimate of $m_j''$ defined in (3.16) involve application of the backfitting regression algorithm. For these pilot estimates one needs to select another set of bandwidths. Iterative schemes to select fully data-dependent plug-in bandwidths are discussed in Section 4.

The next theorem states the conditions under which $\widehat{m}_j''$ is uniformly consistent. This immediately implies that $\widehat{h}_{PL} - h_{\mathrm{opt}}$ and $\widehat{h}_{PL}^* - h_{\mathrm{opt}}^*$ are of lower order.

THEOREM 3.3. *Suppose that assumption* (A5), *in addition to the assumptions of Theorem* 3.1, *holds. Then, for* $g_j$ *with* $g_j \to 0$ *and* $g_j^{-2} n^{-2/5}(\log n)^{1/2} \to 0$, *we have uniformly for* $0 \le u_j \le 1$,

$$\widehat{m}_j''(u_j) - m_j''(u_j) = o_p(1).$$

*Suppose additionally that*

$$\frac{1}{n}\sum_{i=1}^n w(X^i)\{\widetilde{m}(X^i) - m_0 - m_1(X_1^i) - \cdots - m_d(X_d^i)\}^2 = o_P(1).$$

*Then*

$$\widehat{h}_{PL} - h_{\mathrm{opt}} = o_p(n^{-1/5}).$$

*If additionally*

$$\frac{1}{n}\sum_{i=1}^n w_j(X_j^i)\{\widetilde{m}(X^i) - m_0 - m_1(X_1^i) - \cdots - m_d(X_d^i)\}^2 = o_P(1),$$

*then*

$$\widehat{h}_{PL}^* - h_{\mathrm{opt}}^* = o_p(n^{-1/5}).$$

We now give a heuristic discussion of the rates of convergence of $(\widehat{h}_{PL,j} - h_{\mathrm{opt},j})/h_{\mathrm{opt},j}$ and $(\widehat{h}_{PL,j}^* - h_{\mathrm{opt},j}^*)/h_{\mathrm{opt},j}^*$. For simplicity we consider only the latter. Similar arguments may be applied to the former. Note that the rate of the latter coincides with that of

$$n^{-1}\sum_{i=1}^{n} w_j(X_j^i)[\widehat{m}_j''(X_j^i)^2 - m_j''(X_j^i)^2].$$

We now suppose that $\widehat{m}_j''(u_j)$ can be decomposed into $m_j''(u_j) + \text{bias}(u_j) + \text{stoch}(u_j)$, where $\text{bias}(u_j)$ is a bias term and $\text{stoch}(u_j)$ is a mean zero part consisting of local and global averages of $\varepsilon^i$. Under higher-order smoothness conditions one may expect an order of $g_j^2$ for $\text{bias}(u_j)$ and an order of $(ng_j^5)^{-1/2}$ for $\text{stoch}(u_j)$. Now

$$\begin{aligned}
n^{-1}\sum_{i=1}^{n} & w_j(X_j^i)\{\widehat{m}_j''(X_j^i)^2 - m_j''(X_j^i)^2\} \\
&= n^{-1}\sum_{i=1}^{n} w_j(X_j^i)\,\text{bias}(X_j^i)^2 + 2n^{-1}\sum_{i=1}^{n} w_j(X_j^i)\,\text{bias}(X_j^i)\,\text{stoch}(X_j^i) \\
&\quad + n^{-1}\sum_{i=1}^{n} w_j(X_j^i)\,\text{stoch}(X_j^i)^2 + 2n^{-1}\sum_{i=1}^{n} w_j(X_j^i) m_j''(X_j^i)\,\text{bias}(X_j^i) \\
&\quad + 2n^{-1}\sum_{i=1}^{n} w_j(X_j^i) m_j''(X_j^i)\,\text{stoch}(X_j^i).
\end{aligned}$$

By standard reasoning, one may find the following rates of convergence for the five terms on the right-hand side of the above equation: $g_j^4, n^{-1/2}g_j^2, n^{-1}g_j^{-5}, g_j^2, n^{-1/2}$. The maximum of these orders is minimized by $g_j \sim n^{-1/7}$, leading to $(\widehat{h}_{PL,j} - h_j^{\text{opt}})/h_j^{\text{opt}} = O_p(n^{-2/7})$.

The relative rate $O_p(n^{-2/7})$ for the plug-in bandwidth selectors is also achieved by the fully automated bandwidth selector of Opsomer and Ruppert [14], and is identical to the rate of the plug-in rule for the one-dimensional local linear regression estimator of Ruppert, Sheather and Wand [18]. We note here that more sophisticated choices of the constant factor of $n^{-1/7}$ for the bandwidth $g_j$ would yield faster rates such as $n^{-4/13}$ or even $n^{-5/14}$. See [15, 16] or [19].

**4. Practical implementation of the bandwidth selectors.** We suggest use of iterative procedures for approximation of $\widehat{h}_{PLS}$, $\widehat{h}_{PL}$ and $\widehat{h}_{PL}^*$. We note that use of $\widehat{h}_{PL}$ and $\widehat{h}_{PL}^*$ is restricted to local linear smooth backfitting. For $\widehat{h}_{PLS}$ we propose use of the iterative smooth backfitting algorithm based on (2.6) for Nadaraya–Watson smoothing and (3.2) for the local linear fit, and updating of the bandwidth $h_j$ when the $j$th additive component is calculated in the iteration step. This can be done by computing $PLS(h)$ for a finite number of $h_j$'s with $h_1,\dots,h_{j-1},h_{j+1},\dots,h_d$ being held fixed, and then by replacing $h_j$ by the minimizing value of $h_j$. Specifically, we suggest the following procedure:

*Step* 0. Initialize $h_j^{[0]}$ for $j=1,\dots,d$.

*Step* r. Find $h_j^{[r]} = \arg\min_{h_j} PLS(h_1^{[r-1]},\dots,h_{j-1}^{[r-1]},h_j,h_{j+1}^{[r-1]},\dots,h_d^{[r-1]})$ on a grid of $h_j$, for $j=1,\dots,d$.

The computing time for the above iterative procedure to find $\widehat{h}_{PLS}$ is $R\times d\times N\times C$ where $R$ denotes the number of iterations, $N$ is the number of points on the grid of each $h_j$ and $C$ is the time for the evaluation of $PLS$ (or equivalently $RSS$) with a given set of bandwidths. This is much less than the computing time required for the $d$-dimensional grid search, which is $N^d\times C$.

In the implementation of the iterative smooth backfitting algorithm, the estimate $\widehat{m}_j$ could be calculated on a grid of $I_j$. The integrals used in the updating steps of the smooth backfitting can be replaced by the weighted averages over this grid. For the calculation of $PLS(h)$ we need $\widehat{m}_j(X_j^i)$. These values can be approximated by linear interpolation between the neighboring points on the grid. In the simulation study presented in the next section we used a grid of 25 equally spaced points in the interval $I_j=[0,1]$.

Next we discuss how to approximate $\widehat{h}_{PL}$ for the local linear smooth backfitting. We calculate the residuals by use of a backfitting estimate. This means that we replace $n^{-1}\sum_{i=1}^n w(X^i)(\widehat{\varepsilon}^i)^2$ in $\widehat{AASE}$ by $RSS$ as defined in (3.5). Recall that $RSS$ involves the bandwidth $h=(h_1,\dots,h_d)$, and that the local quadratic estimate $\widehat{m}_j''$ defined in (3.16) depends on the bandwidth $g_j$ as well as $h$. The residual sum of squares $RSS$ and the estimate $\widehat{m}_j''$ depend on the bandwidths $h$ of the smooth backfitting regression estimates. To stress this dependence on $h$ and $g_j$, we write $RSS(h)$ and $\widehat{m}_j''(\cdot;h,g_j)$ for $RSS$ and $\widehat{m}_j''$, respectively. We propose the following iterative procedure for $\widehat{h}_{PL}$:

*Step* 0. Initialize $h^{[0]}=(h_1^{[0]},\dots,h_d^{[0]})$.

*Step* r. Compute on a grid of $h=(h_1,\dots,h_d)$

$$\begin{aligned}\widehat{AASE}^{[r]}(h) &= RSS(h^{[r-1]})\left\{\int K^2(t)\,dt\right\}\sum_{j=1}^d \frac{1}{nh_j}\\ &\quad + \frac{1}{4n}\sum_{i=1}^n w(X^i)\left\{\sum_{j=1}^d h_j^2\widehat{m}_j''(X_j^i;h^{[r-1]},g_j^{[r-1]})\right\}^2\\ &\quad \times\left\{\int t^2K(t)\,dt\right\}^2\end{aligned}$$

with $g_j^{[r-1]} = c\,h_j^{[r-1]}$ ($c = 1.5$ or 2, say), and then find

$$h^{[r]} = \arg\min \widehat{AASE}^{[r]}(h).$$

A more sophisticated choice of $g_j$ suggested by the discussion at the end of Section 3 would be $g_j = ch_j^{5/7}$ for some properly chosen constant $c > 0$.

We also give an alternative algorithm to approximate $\widehat{h}_{PL}$, which requires only a one-dimensional grid search. This would be useful for very high-dimensional covariates:

*Step* 0′. Initialize $h^{[0]} = (h_1^{[0]}, \dots, h_d^{[0]})$.
*Step* r′. For $j = 1, \dots, d$, compute on a grid of $h_j$

$$\begin{aligned}
&\widehat{AASE}^{[r,j]}(h_1^{[r-1]}, \dots, h_{j-1}^{[r-1]}, h_j, h_{j+1}^{[r-1]}, \dots, h_d^{[r-1]}) \\
&\quad = RSS(h^{[r-1]})\left\{\int K^2(t)\,dt\right\}\left\{\frac{1}{nh_j} + \sum_{\ell \neq j}^{d} \frac{1}{nh_\ell^{[r-1]}}\right\} \\
&\qquad + \frac{1}{4n}\sum_{i=1}^{n} w(X^i)\left\{h_j^2 \widehat{m}_j''(X_j^i; h^{[r-1]}, g_j^{[r-1]})\right. \\
&\qquad\qquad \left. + \sum_{\ell \neq j}^{d} (h_\ell^{[r-1]})^2 \widehat{m}_\ell''(X_\ell^i; h^{[r-1]}, g_\ell^{[r-1]})\right\}^2 \\
&\qquad \times \left\{\int t^2 K(t)\,dt\right\}^2,
\end{aligned}$$

and then find

$$h_j^{[r]} = \arg\min_{h_j} \widehat{AASE}^{[r,j]}(h_1^{[r-1]}, \dots, h_{j-1}^{[r-1]}, h_j, h_{j+1}^{[r-1]}, \dots, h_d^{[r-1]}).$$

In the grid search for $h_j^{[r]}$ we use $h^{[r-1]}$ rather than $(h_1^{[r-1]}, \dots, h_{j-1}^{[r-1]}, h_j, h_{j+1}^{[r-1]}, \dots, h_d^{[r-1]})$ for $RSS$ and $\widehat{m}_\ell''$. The reason is that the latter requires repetition of the whole backfitting procedure (3.2) for every point on the grid of the bandwidth. Thus, it is computationally much more expensive than our second proposal for approximating $\widehat{h}_{PL}$.

Finally, we give an algorithm to approximate $\widehat{h}_{PL,j}^*$. We suppose calculation of the residuals by use of a backfitting estimate. This means that we replace $n^{-1}\sum_{i=1}^{n} w_j(X_j^i)(\widehat{\varepsilon}^i)^2$ in (3.15) by $RSS$. Thus, $\widehat{h}_{PL,j}^*$ is given by

$$\begin{aligned}
\widehat{h}_{PL,j}^* &= n^{-1/5}\left[RSS \times \left\{\int K^2(t)\,dt\right\}\right]^{1/5} \\
&\quad \times \left[n^{-1}\sum_{i=1}^{n} w_j(X_j^i)\widehat{m}_j''(X_j^i)^2 \left\{\int t^2 K(t)\,dt\right\}^2\right]^{-1/5}.
\end{aligned} \tag{4.1}$$

We propose the following iterative procedure for $\widehat{h}^*_{PL}$. Start with some initial bandwidths $h_1, \ldots, h_d$ and calculate $\widehat{m}_1^{LL}, \ldots, \widehat{m}_d^{LL}$ with these bandwidths, and compute $RSS$. Choose $g_j = ch_j$ (with $c = 1.5$ or 2, say). Then calculate $\widehat{m}''_1, \ldots, \widehat{m}''_d$ by (3.16). Plug $RSS$ and the computed values of $\widehat{m}''_j(X^i_j)$'s into (4.1), which defines new bandwidths $h_1, \ldots, h_d$. Then the procedure can be iterated.

It was observed in the simulation study presented in Section 5 that the iterative algorithms for approximating $\widehat{h}_{PLS}$, $\widehat{h}_{PL}$ and $\widehat{h}^*_{PL}$ converge very quickly. With the convergence criterion $10^{-3}$ on the relative changes of the bandwidth selectors, the average (out of 500 cases) numbers of iterations for the three bandwidth selectors were 4.27, 6.30 and 5.23, respectively. The worst cases had eight iterations.

**5. Simulations.** In this section we present simulations for the small sample performance of the plug-in selectors $\widehat{h}_{PL}$, $\widehat{h}^*_{PL}$ and the penalized least squares bandwidth $\widehat{h}_{PLS}$. We will do this only for local linear smooth backfitting.

Our first goal was to compare how much these bandwidths differ from their theoretical targets. For this, we drew 500 datasets $(X^i, Y^i)$, $i = 1, \ldots, n$, with $n = 200$ and 500 from the model

$$Y^i = m_1(X^i_1) + m_2(X^i_2) + m_3(X^i_3) + \varepsilon^i, \qquad (M1)$$

where $m_1(x_1) = x_1^2$, $m_2(x_2) = x_2^3$, $m_3(x_3) = x_3^4$ and $\varepsilon^i$ are distributed as $N(0, 0.01)$. The covariate vectors were generated from joint normal distributions with marginals $N(0.5, 0.5)$ and correlations $(\rho_{12}, \rho_{13}, \rho_{23}) = (0, 0, 0)$ and $(0.5, 0.5, 0.5)$. Here $\rho_{ij}$ denotes the correlation between $X_i$ and $X_j$. If the generated covariate vector was within the cube $[0, 1]^3$, then it was retained in the sample. Otherwise, it was removed. This was done until arriving at the predetermined sample size 200 or 500. Thus, the covariate vectors follow truncated normal distributions and have compact support $[0, 1]^3$ satisfying assumption (A2). Both of the kernels $K$ that we used for the backfitting algorithm and $L$ for estimating $m''_j$ by (3.16) were the biweight kernel $K(u) = L(u) = (15/16)(1 - u^2)^2 I_{[-1,1]}(u)$. The weight function $w$ in (3.5) and (3.6) was the indicator $\mathbf{1}(u \in [0, 1]^3)$, and $w_j$ in (3.13) and (4.1) was $\mathbf{1}(u_j \in [0, 1])$.

Kernel density estimates of the densities of $\log(\widehat{h}_{PLS,j}) - \log(\widehat{h}_{ASE,j})$, $\log(\widehat{h}_{PL,j}) - \log(\widehat{h}_{ASE,j})$ and $\log(\widehat{h}^*_{PL,j}) - \log(\widehat{h}_{ASE,j})$ are overlaid in Figures 1–3 for $j = 1, 2, 3$. The results are based on 500 replicates for the two choices of the correlation values and of the sample sizes. The kernel density estimates were constructed by using the standard normal kernel and the common bandwidth 0.12. The iterative procedures described in Section 4 for $\widehat{h}_{PLS}$, $\widehat{h}_{PL}$ and $\widehat{h}^*_{PL}$ were used here. In all cases, the initial bandwidth $h^{[0]} = (0.1, 0.1, 0.1)$ was used. For $\widehat{h}_{PL}$, the first proposal with three-dimensional grid search was implemented. We tried $g = 1.5h$ and $g = 2h$ to estimate $m''_j$ in the iterative procedures. We found there is little difference between these two choices, and thus present here only the results for the case $g = 1.5h$. In each of Figures 1–3, the upper two panels show the densities of the log differences for the sample size $n = 200$, while the lower two correspond to the cases where $n = 500$.

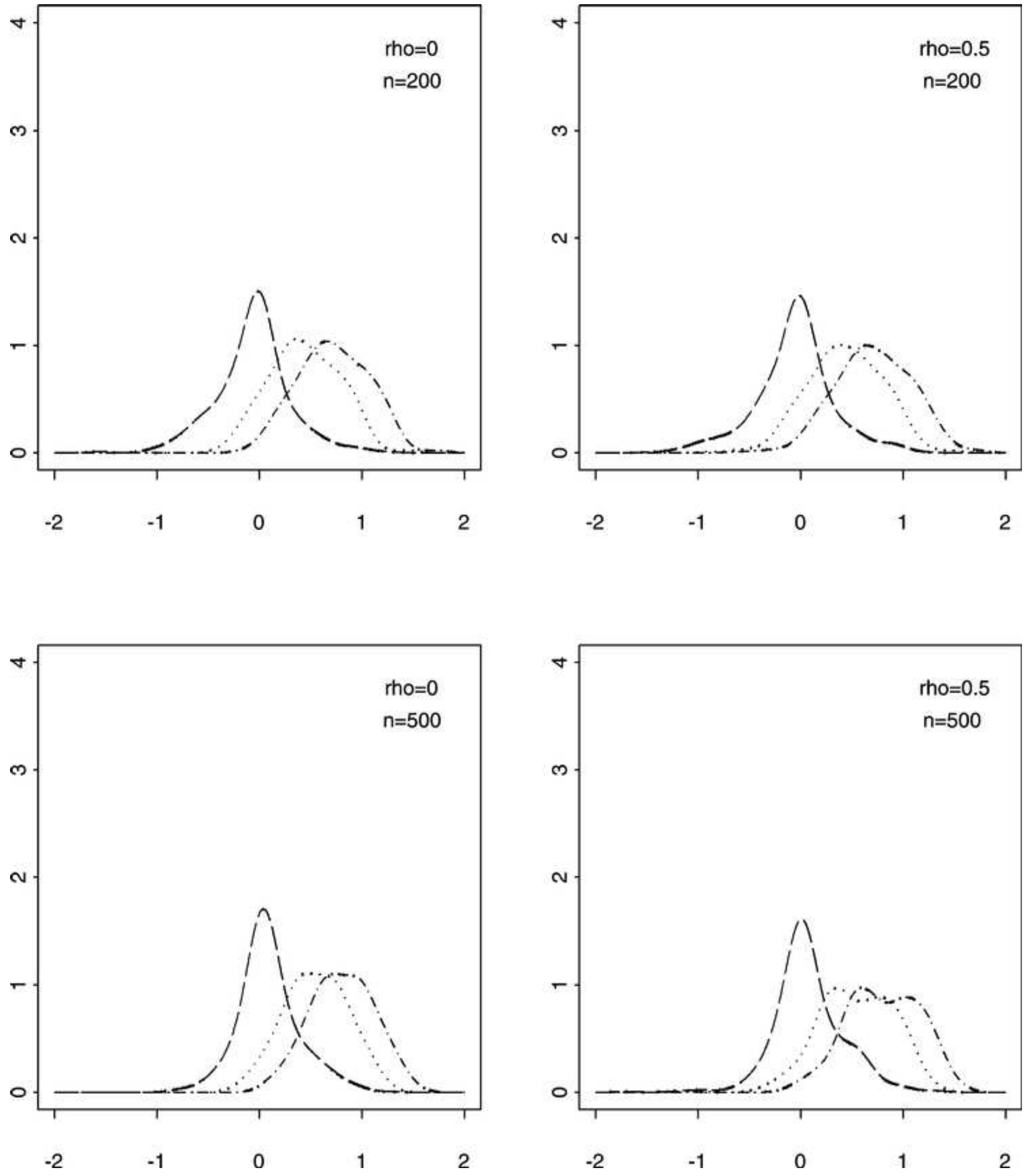

FIG. 1. *Densities of* $\log(\widehat{h}_1) - \log(\widehat{h}_{ASE,1})$ *constructed by the kernel method based on* 500 *pseudosamples. The long-dashed, dotted and dot-dashed curves correspond to* $\widehat{h}_1 = \widehat{h}_{PLS,1}$, $\widehat{h}_{PL,1}$ *and* $\widehat{h}^*_{PL,1}$, *respectively.*

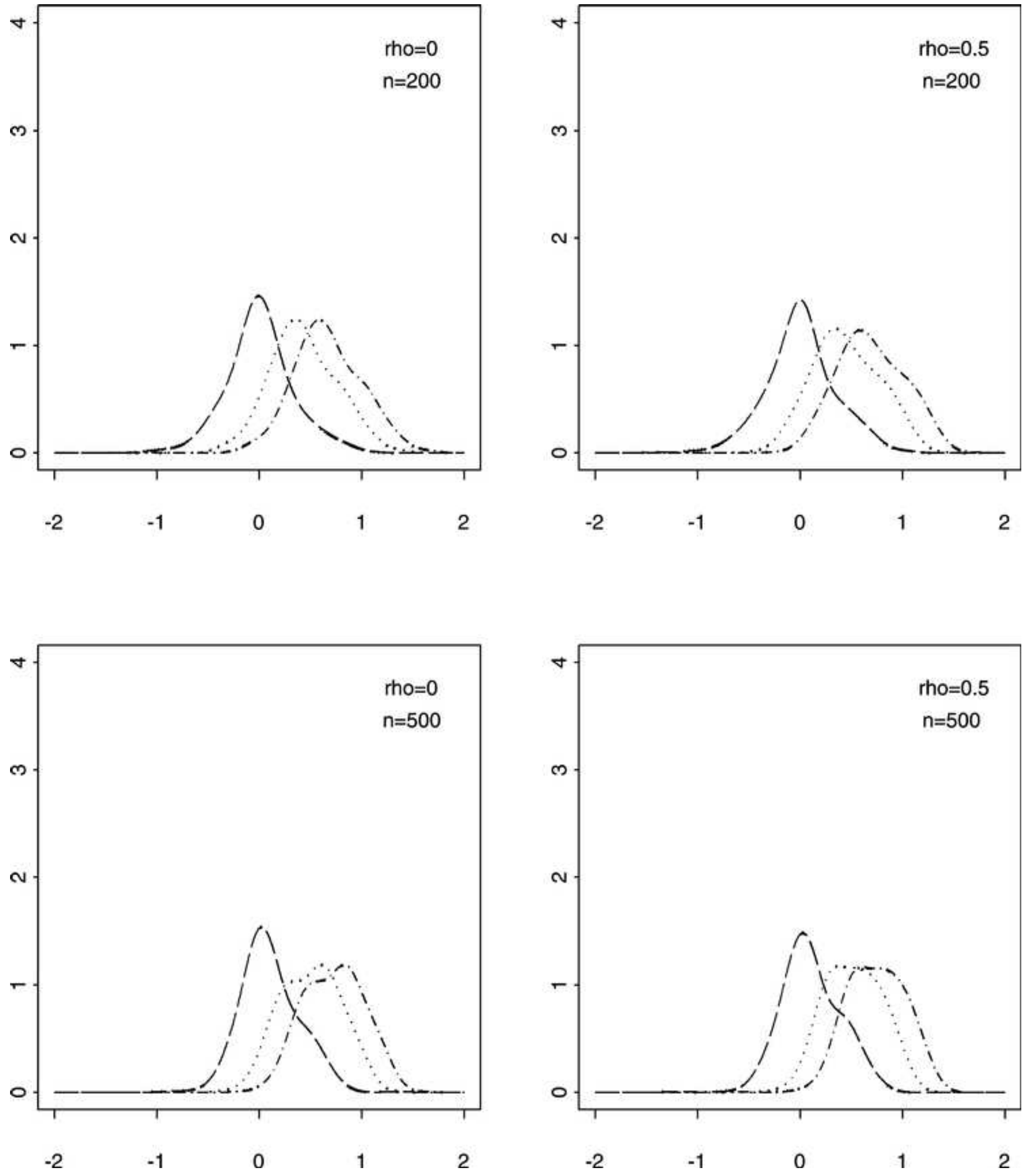

FIG. 2. *Densities of* $\log(\widehat{h}_2) - \log(\widehat{h}_{ASE,2})$ *constructed by the kernel method based on* 500 *pseudosamples. The long-dashed, dotted and dot-dashed curves correspond to* $\widehat{h}_{PLS,2}$, $\widehat{h}_{PL,2}$ *and* $\widehat{h}^*_{PL,2}$, *respectively.*

Comparing the three bandwidth selectors $\widehat{h}_{PLS}$, $\widehat{h}_{PL}$ and $\widehat{h}^*_{PL}$, one sees that the penalized least squares bandwidth has the correct center while the two plug-in bandwidths are positively biased toward $\widehat{h}_{ASE}$. Furthermore, $\widehat{h}_{PLS}$ are less variable than $\widehat{h}_{PL}$ and $\widehat{h}^*_{PL}$ as an estimator of $\widehat{h}_{ASE}$. This shows the penalized least squares approach is superior to the other two methods in terms of estimating $\widehat{h}_{ASE}$. We found, however, $\widehat{h}_{PL}$ and $\widehat{h}^*_{PL}$ are more stable and less biased as estimators of $h_{\text{opt}}$ and $h^*_{\text{opt}}$, respectively.

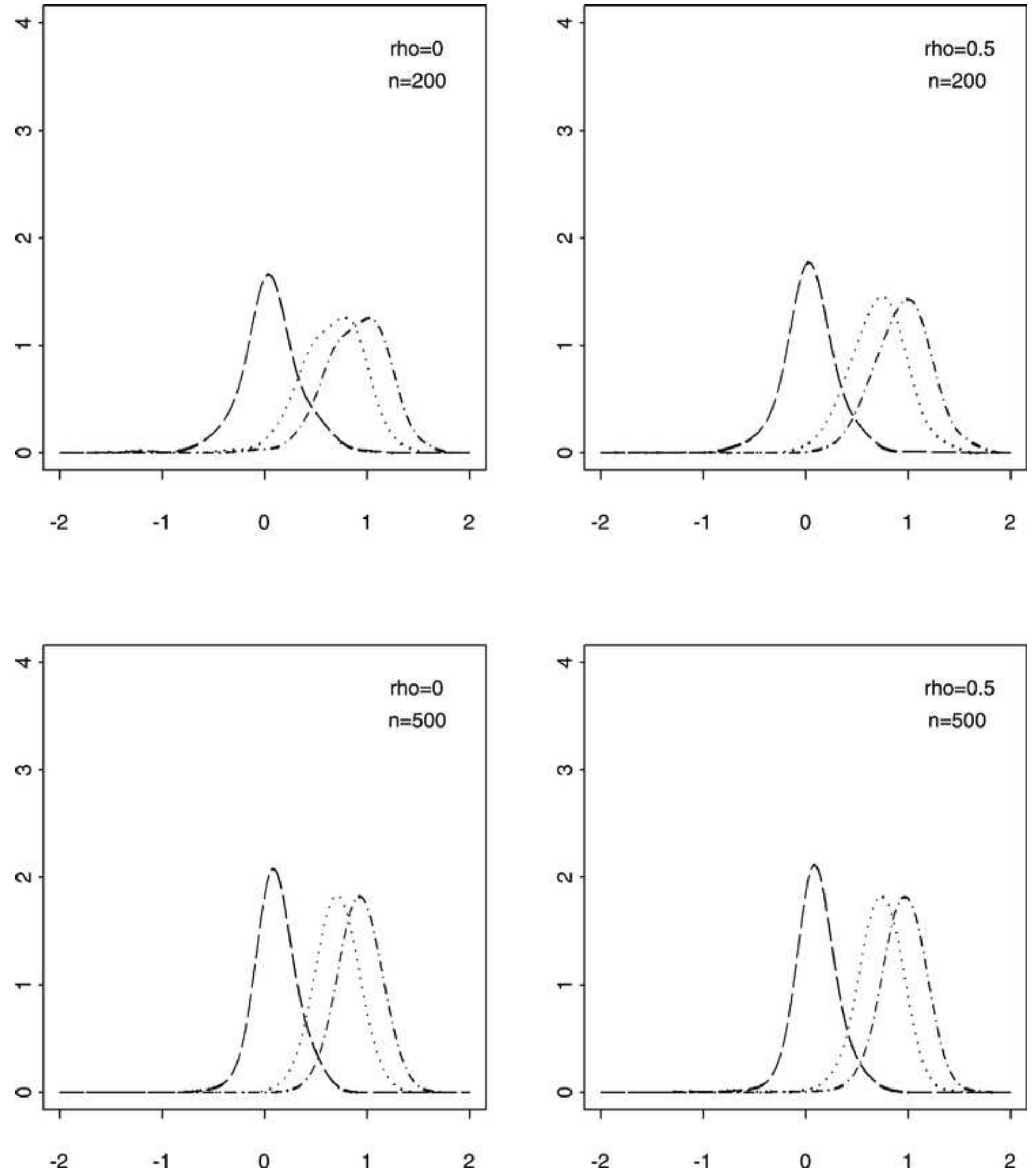

Fig. 3. *Densities of* $\log(\widehat{h}_3^*) - \log(\widehat{h}_{ASE,3})$ *constructed by the kernel method based on* 500 *pseudosamples. The long-dashed, dotted and dot-dashed curves correspond to* $\widehat{h}_{PLS,3}$, $\widehat{h}_{PL,3}$ *and* $\widehat{h}^*_{PL,3}$, *respectively.*

Table 1
*Averages of $ASE(\widehat{h})$ and $ASE_j(\widehat{h})$ $(j=1,2,3)$ for $\widehat{h}=\widehat{h}_{PLS}$, $\widehat{h}_{PL}$ and $\widehat{h}^*_{PL}$, based on 500 pseudosamples*

| | | | $\widehat{h}_{PLS}$ | $\widehat{h}_{PL}$ | | $\widehat{h}^*_{PL}$ | |
|---|---|---|---|---|---|---|---|
| | | | | $g=1.5h$ | $g=2h$ | $g=1.5h$ | $g=2h$ |
| Average | $n=200$ | $\rho=0$ | 0.00251 | 0.00347 | 0.00350 | 0.00471 | 0.00478 |
| $ASE$ | | $\rho=0.5$ | 0.00247 | 0.00362 | 0.00367 | 0.00513 | 0.00521 |
| | $n=500$ | $\rho=0$ | 0.00130 | 0.00195 | 0.00199 | 0.00269 | 0.00277 |
| | | $\rho=0.5$ | 0.00133 | 0.00209 | 0.00213 | 0.00294 | 0.00303 |
| Average | $n=200$ | $\rho=0$ | 0.00107 | 0.00131 | 0.00133 | 0.00169 | 0.00172 |
| $ASE_1$ | | $\rho=0.5$ | 0.00112 | 0.00150 | 0.00153 | 0.00207 | 0.00211 |
| | $n=500$ | $\rho=0$ | 0.00045 | 0.00063 | 0.00065 | 0.00084 | 0.00088 |
| | | $\rho=0.5$ | 0.00052 | 0.00076 | 0.00078 | 0.00103 | 0.00108 |
| Average | $n=200$ | $\rho=0$ | 0.00104 | 0.00085 | 0.00085 | 0.00078 | 0.00078 |
| $ASE_2$ | | $\rho=0.5$ | 0.00100 | 0.00079 | 0.00079 | 0.00072 | 0.00072 |
| | $n=500$ | $\rho=0$ | 0.00044 | 0.00037 | 0.00037 | 0.00033 | 0.00033 |
| | | $\rho=0.5$ | 0.00047 | 0.00038 | 0.00038 | 0.00034 | 0.00034 |
| Average | $n=200$ | $\rho=0$ | 0.00112 | 0.00079 | 0.00079 | 0.00073 | 0.00073 |
| $ASE_3$ | | $\rho=0.5$ | 0.00121 | 0.00090 | 0.00090 | 0.00086 | 0.00086 |
| | $n=500$ | $\rho=0$ | 0.00051 | 0.00038 | 0.00037 | 0.00034 | 0.00033 |
| | | $\rho=0.5$ | 0.00061 | 0.00050 | 0.00050 | 0.00047 | 0.00047 |

It is also interesting to compare the performance of the bandwidth selectors in terms of the average squared error of the resulting regression estimator. Table 1 shows the means (out of 500 cases) of the $ASE$ and $ASE_j$ for the three bandwidth selectors. First, it is observed that $\widehat{h}_{PLS}$ produces the least $ASE$. This means that $\widehat{h}_{PLS}$ is most effective for estimating the whole regression function. Now, for accuracy of each one-dimensional component estimator, none of the bandwidth selectors dominates the others in all cases. For $ASE_1$, the penalized least squares bandwidth does the best, while for $ASE_2$ and $ASE_3$ the plug-in $\widehat{h}^*_{PL}$ shows the best performance. The backfitting estimates the centered true component functions because of the normalization (3.4). Thus, $\widehat{m}_j^{LL}(x_j)$ estimates $m_j(x_j) - Em_j(X_j^i)$, not $m_j(x_j)$. We used these centered true functions to compute $ASE_j$.

Table 1 also shows that the means of the average squared errors are reduced approximately by half when the sample size is increased from 200 to 500. Although not reported in the table, we computed $E(\widehat{h}_{j,200})/E(\widehat{h}_{j,500})$ for the three bandwidth selectors, where $\widehat{h}_{j,n}$ denotes the bandwidth selector for the $j$th component from a sample of size $n$. We found that these values vary within the range $(1.20, 1.26)$ which is roughly $(200/500)^{-1/5}$. This means the assumed rate $O(n^{-1/5})$ for the bandwidth selectors actually holds in practice. Now, we note that the increase of correlation from 0 to 0.5 does not deteriorate much the means of the $ASE$ and $ASE_j$. However, we found in a separate experiment that in a more extremal case of $\rho_{ij} \equiv 0.9$ the means of the $ASE$ and $ASE_j$ are increased by a factor of 3 or 4. In another separate experiment where the noise level is 0.1, that is, the errors are generated from $N(0, 0.1)$, we observed that the means of the $ASE$ and $ASE_j$ are increased by a factor of 3 or 4, too. The main lessons on comparison of the three bandwidth selectors from these two separate experiments are the same as in the previous paragraph.

Figure 4 visualizes the overall performance of the backfitting for the three bandwidth selectors. For each $\widehat{h} = \widehat{h}_{ASE}$, $\widehat{h}_{PLS}$, $\widehat{h}_{PL}$, $\widehat{h}^*_{PL}$, we computed $ASE(\widehat{h})$ and $ASE_j(\widehat{h})$ for 500 datasets and arranged the 500 values of $d = ASE(\widehat{h})$ or $ASE_j(\widehat{h})$ in increasing order. Call them $d_{(1)} \leq d_{(2)} \leq \cdots \leq d_{(500)}$. Figure 4 shows the quantile plots $\{i/500, d_{(i)}\}_{i=1}^{500}$ for the case where $n = 500$ and $\rho_{ij} \equiv 0.5$. The bandwidth $g = 1.5h$ was used in the pilot estimation step for the two plug-in bandwidths. The figure reveals that the quantile function of $ASE(\widehat{h})$ for $\widehat{h} = \widehat{h}_{PLS}$ is consistently below those for the two plug-in rules and is very close to that for $\widehat{h} = \widehat{h}_{ASE}$. For $ASE_j(\widehat{h})$, none of the three bandwidth selectors dominates the other two for all $j$, the result also seen in Table 1, but in any case the quantile function of $ASE_j(\widehat{h}_{PLS})$ is closest to that of $ASE_j(\widehat{h}_{ASE})$. We note that the quantile functions of $ASE_j(\widehat{h}_{ASE})$ are not always the lowest since $\widehat{h}_{ASE} = (\widehat{h}_{ASE,1}, \widehat{h}_{ASE,2}, \widehat{h}_{ASE,3})$ does not minimize each component's $ASE_j$.

Asymptotic theory says that in first order the accuracy of the backfitting estimate does not decrease with increasing number of additive components. And this also holds for the backfitting estimates with the data-adaptively chosen bandwidths. We wanted to check if this also holds for finite samples. For this purpose we compared our model ($M1$) with three additive components with a model that has only one component,

$$(M2) \qquad Y^i = m_1(X_1^i) + \varepsilon^i.$$

We drew 500 datasets of sizes 200 and 500 from the models ($M1$) and ($M2$). The errors and the covariates at the correlation level 0.5 were generated in the same way as described in the second paragraph of this section. The penalized least squares bandwidth for the single covariate case, denoted by $\widehat{h}_{PLS(1)}$, was obtained by minimizing

$$PLS_1(h_1) = RSS_1(h_1)\left\{1 + \frac{2K(0)}{nh_1}\right\},$$

where $RSS_1(h_1) = n^{-1}\sum_{i=1}^n \{Y^i - \widetilde{m}_1^{LL}(X_1^i; h_1)\}^2$ and $\widetilde{m}_1^{LL}(\cdot; h_1)$ is the ordinary local linear fit for model ($M2$) with bandwidth $h_1$. The plug-in bandwidth selector, $\widehat{h}_{PL(1)}$, for the single covariate case was obtained by a formula

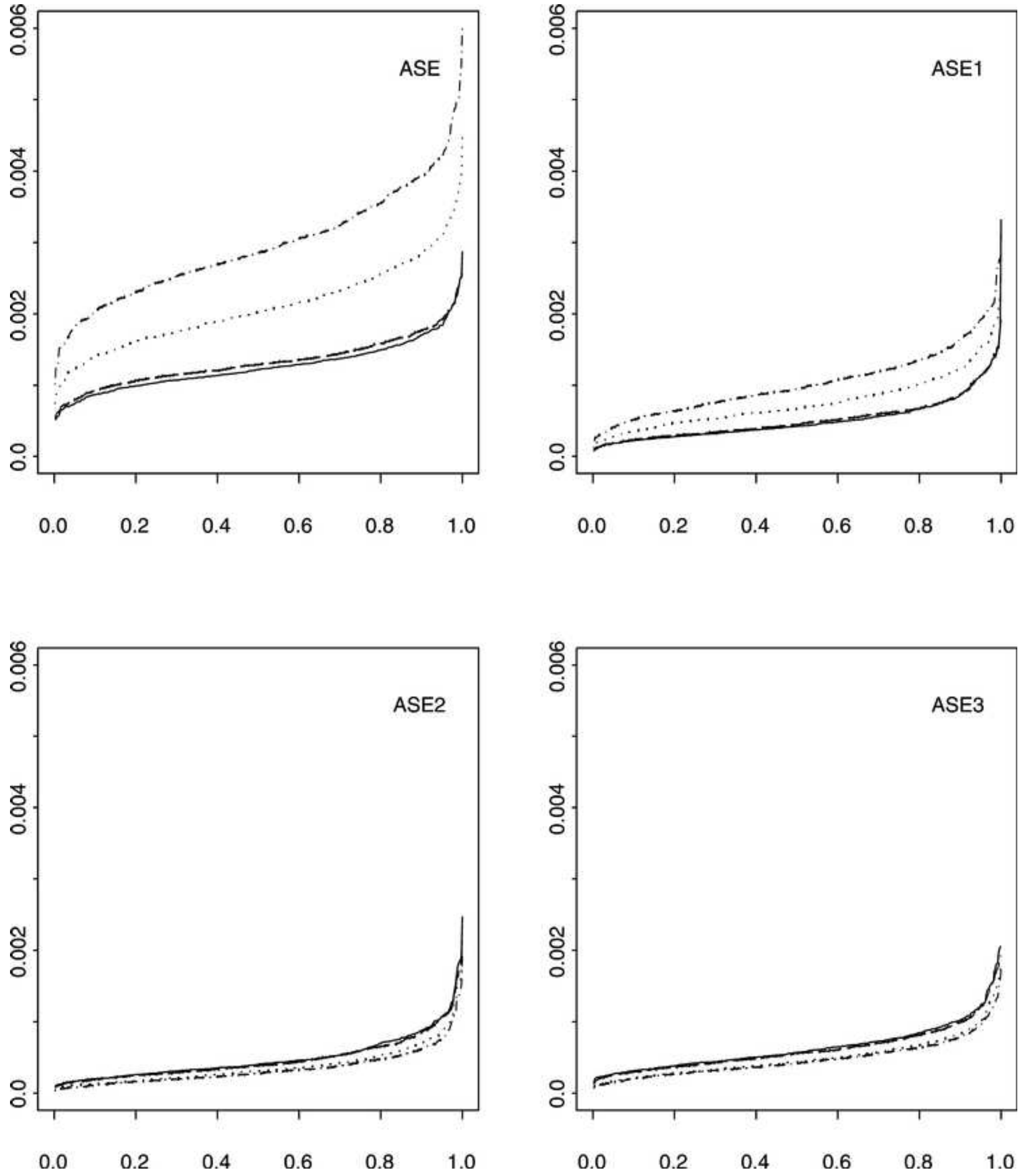

FIG. 4. *Quantile functions of* $ASE(\widehat{h})$ *and* $ASE_j(\widehat{h})$ *for* $\widehat{h}_{ASE}$ *and the three bandwidth selectors. Solid, long-dashed, dotted and dot-dashed curves correspond to* $\widehat{h}_{ASE}$, $\widehat{h}_{PLS}$, $\widehat{h}_{PL}$ *and* $\widehat{h}^*_{PL}$, *respectively. The sample size was* $n = 500$ *and the correlations between the covariates were all* 0.5.

similar to the one in (4.1), where $RSS$ is replaced by $RSS_1$ and $\widetilde{m}_1^{LL}$, instead of the backfitting estimate $\widehat{m}_1^{LL}$, is used to calculate the local quadratic estimate of $m_1''$. For $\widehat{h}_{PL(1)}$, an iterative procedure similar to those described in Section 4 was used here, again with the choice $g = 1.5h$. Table 2 shows $E\{ASE_1(\widehat{h})\}$ for $\widehat{h} = \widehat{h}_{PLS}$, $\widehat{h}_{PLS(1)}$, $\widehat{h}_{PL}$, $\widehat{h}^*_{PL}$ and $\widehat{h}_{PL(1)}$. Also, it gives the

TABLE 2

*Averages of $ASE_1(\widehat{h})$ as an error criterion for estimating the first component $m_1$, based on 500 pseudosamples from the models (M1) and (M2)*

| | | $\hat{h}_{PLS}$ | $\hat{h}_{PLS(1)}$ | $\hat{h}_{PL}$ | $\hat{h}^*_{PL}$ | $\hat{h}_{PL(1)}$ |
|---|---|---|---|---|---|---|
| $n = 200$ | $\rho = 0$ | 0.00107 (3.147) | 0.00034 | 0.00131 (4.517) | 0.00169 (5.828) | 0.00029 |
| | $\rho = 0.5$ | 0.00112 (3.394) | 0.00033 | 0.00150 (5.357) | 0.00207 (7.393) | 0.00028 |
| $n = 500$ | $\rho = 0$ | 0.00045 (3.000) | 0.00015 | 0.00063 (4.500) | 0.00084 (6.000) | 0.00014 |
| | $\rho = 0.5$ | 0.00052 (3.714) | 0.00014 | 0.00076 (5.846) | 0.00103 (7.923) | 0.00013 |

Also given in the parentheses are the relative increases of $E\{ASE_1(\widehat{h})\}$ due to the increased dimension of the covariates. The choice $g = 1.5h$ was used for the plug-in rules.

relative increases of $E\{ASE_1(\widehat{h})\}$ due to the increased dimension of the covariates. For $\widehat{h}_{PLS(1)}$ and $\widehat{h}_{PL(1)}$ the one-dimensional local linear estimate $\widetilde{m}_1^{LL}$ and the noncentered regression function $m_1$ were used to compute the values of $ASE_1$.

From Table 2, it appears that the increased dimension of the covariates has some considerable effect on the regression estimates. The relative increase of $ASE_1$ for the penalized least squares bandwidth is smaller than those for the plug-in rules, however. Also, one observes higher rates of increase for the correlated covariates. An interesting fact is that $\widehat{h}_{PL(1)}$ is slightly better than $\widehat{h}_{PLS(1)}$ in the single covariate case. The results for the other component functions, which are not presented here, showed the same qualitative pictures.

## 6. Assumptions, auxiliary results and proofs.

6.1. *Assumptions.* We use the following assumptions.

(A1) The kernel $K$ is bounded, has compact support ($[-1.1]$, say), is symmetric about zero and is Lipschitz continuous, that is, there exists a positive finite constant $C$ such that $|K(t_1) - K(t_2)| \le C|t_1 - t_2|$.

(A2) The $d$-dimensional vector $X^i$ has compact support $I_1 \times \cdots \times I_d$ for bounded intervals $I_j$, and its density $p$ is bounded away from zero and infinity on $I_1 \times \cdots \times I_d$. The tuples $(X^i, \varepsilon^i)$ are i.i.d.

(A3) Given $X^i$ the error variable $\varepsilon^i$ has conditional zero mean, and for some $\gamma > 4$ and $C' < \infty$

$$E[|\varepsilon^i|^\gamma | X^i] < C' \qquad \text{a.s.}$$

(A4) The functions $m_j'', p_j'$ and $(\partial/\partial x_j) p_{jk}(x_j, x_k)$ $(1 \le j, k \le d)$ exist and are continuous.

(A5) The kernel $L$ is twice continuously differentiable and has bounded support ($[-1, 1]$, say).

6.2. *Auxiliary results.* In this section we will give higher-order expansions of $\widehat{m}_j^{NW}$ and $\widehat{m}_j^{LL}$. These expansions will be used in the proofs of Theorems 2.1, 2.2 and 3.3. The expansions given in [9] are only of order $o_p(n^{-2/5})$. Furthermore, they are not uniform in $h$. For the proof of our results we need expansions of order $o_p(n^{-1/2})$. First, we consider the Nadaraya–Watson smooth backfitting estimate $\widehat{m}_j^{NW}$.

As in [9] we decompose $\widehat{m}_j^{NW}$ into

$$\widehat{m}_j^{NW}(u_j) = \widehat{m}_j^{NW,A}(u_j) + \widehat{m}_j^{NW,B}(u_j),$$

where $\widehat{m}_j^{NW,S}$ $(S = A, B)$ is defined by

$$\begin{aligned} \widehat{m}_j^{NW,S}(u_j) &= \widetilde{m}_j^{NW,S}(u_j) \\ &\quad - \sum_{k \ne j} \int_{I_k} \widehat{m}_k^{NW,S}(u_k) \frac{\widehat{p}_{j,k}(u_j, u_k)}{\widehat{p}_j(u_j)} \, du_k - \widehat{m}_0^{NW,S}, \end{aligned} \tag{6.1}$$

where $\widehat{m}_0^{NW,A} = n^{-1} \sum_{i=1}^n \varepsilon^i$, $\widehat{m}_0^{NW,B} = n^{-1} \sum_{i=1}^n \{m_0 + \sum_{j=1}^d m_j(X_j^i)\}$ and

$$\widetilde{m}_j^{NW,A}(u_j) = \widehat{p}_j(u_j)^{-1} n^{-1} \sum_{i=1}^n K_{h_j}(u_j, X_j^i) \varepsilon^i,$$

$$\widetilde{m}_j^{NW,B}(u_j) = \widehat{p}_j(u_j)^{-1} n^{-1} \sum_{i=1}^n K_{h_j}(u_j, X_j^i) \left\{ m_0 + \sum_{j=1}^d m_j(X_j^i) \right\}.$$

Here $\widetilde{m}_j^{NW,B}$ and $\widehat{m}_j^{NW,B}$ are related to the sum of the true function and the bias, whereas $\widetilde{m}_j^{NW,A}$ and $\widehat{m}_j^{NW,A}$ represent the "stochastic" part. In particular, $\widetilde{m}_j^{NW,B}$ and $\widehat{m}_j^{NW,B}$ do not depend on the error variables.

We now state our stochastic expansions of $\widehat{m}_j^{NW,A}$ and $\widehat{m}_j^{NW,B}$.

THEOREM 6.1. *Suppose that the assumptions of Theorem 2.1 apply, and that $\widehat{m}_j^{NW,A}$ and $\widehat{m}_j^{NW,B}$ are defined according to (6.1). Then there exist random variables $R_{n,i,j}(u_j, h, X)$, depending on $0 \le u_j \le 1$, $h = (h_1, \ldots, h_d)$ and $X = (X^1, \ldots, X^n)$ (but not on $\varepsilon$), such that*

$$\widehat{m}_j^{NW,A}(u_j) = \widetilde{m}_j^{NW,A}(u_j) + n^{-1} \sum_{i=1}^n R_{n,i,j}(u_j, h, X) \varepsilon^i, \tag{6.2a}$$

$$\text{(6.2b)} \qquad \sup_{0\le u_j\le 1}\ \sup_{C_H n^{-1/5}\le h_1,\dots,h_d\le C_H' n^{-1/5}} |R_{n,i,j}(u_j,h,X)| = O_p(1),$$

$$\text{(6.2c)} \qquad \begin{aligned} &\sup_{0\le u_j\le 1}\ \sup_{C_H n^{-1/5}\le h_1,h_1',\dots,h_d,h_d'\le C_H' n^{-1/5}} |R_{n,i,j}(u_j,h',X) - R_{n,i,j}(u_j,h,X)| \\ &\quad = \sum_{j=1}^d |h_j' - h_j| O_p(n^\alpha) \qquad \text{for some } \alpha > 0. \end{aligned}$$

*Furthermore, uniformly for* $C_H n^{-1/5} \le h_1,\dots,h_d \le C_H' n^{-1/5}$ *and* $0 \le u_j \le 1$,

$$\text{(6.3)} \qquad \widehat{m}_j^{NW,A}(u_j) = \widetilde{m}_j^{NW,A}(u_j) + n^{-1}\sum_{i=1}^n r_{ij}(u_j)\varepsilon^i + o_p(n^{-1/2}),$$

$$\text{(6.4)} \qquad \widehat{m}_j^{NW,B}(u_j) = m_j(u_j) + O_p(n^{-1/5}),$$

*where* $r_{ij}$ *are absolutely uniformly bounded functions with*

$$\text{(6.5)} \qquad |r_{ij}(u_j') - r_{ij}(u_j)| \le C|u_j' - u_j|$$

*for a constant* $C > 0$. *In particular, uniformly for* $C_H n^{-1/5} \le h_j \le C_H' n^{-1/5}$ *and* $h_j \le u_j \le 1 - h_j$,

$$\text{(6.6)} \qquad \widehat{m}_j^{NW,B}(u_j) = m_j(u_j) + \beta_j(h,u_j) + o_p(n^{-2/5}),$$

*where* $\beta_j$ *is chosen so that*

$$\begin{aligned} &\int \beta_j(h,u_j)p_j(u_j)\,du_j \\ &\quad = -\gamma_{n,j} \\ &\quad = \tfrac{1}{2}h_j^2 \int [m_j'(x_j)p_j'(x_j) + \tfrac{1}{2}m_j''(x_j)p_j(x_j)]\,dx_j \int u^2K(u)\,du. \end{aligned}$$

*This choice is possible because of* $\int \beta(h,x)p(x)\,dx = -\sum_{j=1}^d \gamma_{n,j}$.

We now come to the local linear smooth backfitting estimate. For a theoretical discussion, we now decompose this estimate into a stochastic and a deterministic term. For $S = A, B$, define $\widehat{m}_j^{LL,S}$ by

$$\text{(6.7)} \qquad \begin{aligned} \begin{pmatrix} \widehat{m}_j^{LL,S}(u_j) \\ \widehat{m}_j^{LL,1,S}(u_j) \end{pmatrix} &= -\begin{pmatrix} \widehat{m}_0^{LL,S} \\ 0 \end{pmatrix} + \begin{pmatrix} \widetilde{m}_j^{LL,S}(u_j) \\ \widetilde{m}_j^{LL,1,S}(u_j) \end{pmatrix} \\ &\qquad - \widehat{M}_j(u_j)^{-1} \sum_{l\ne j} \int \widehat{S}_{lj}(u_l,u_j) \begin{pmatrix} \widehat{m}_l^{LL,S}(u_l) \\ \widehat{m}_l^{LL,1,S}(u_l) \end{pmatrix} du_l, \end{aligned}$$

$$\widehat{m}_0^{LL,S} = n^{-1} \sum_{i=1}^{n} Y^{i,S}, \tag{6.8}$$

$$\int \widehat{m}_j^{LL,S}(u_j)\widehat{p}_j(u_j)\,du_j + \int \widehat{m}_j^{LL,1,S}(u_j)\widehat{p}_j^1(u_j)\,du_j = 0, \qquad j = 1,\ldots,d, \tag{6.9}$$

where $Y^{i,S} = \varepsilon^i$ for $S = A$ and $m_0 + \sum_{j=1}^d m_j(X_j^i)$ for $S = B$. Furthermore, $\widetilde{m}_j^{LL,S}$ and $\widetilde{m}_j^{LL,1,S}$ are the local linear estimates of the function itself and its first derivative, respectively, for the regression of $\varepsilon^i$ (for $S = A$) or $m_0 + m_1(X_1^i) + \cdots + m_d(X_d)^i$ (for $S = B$) onto $X_j^i$.

For the local linear smooth backfitting estimate, we get the following stochastic expansions.

THEOREM 6.2. *Suppose that the assumptions of Theorem 3.1 apply, and that $\widehat{m}_j^{LL,S}$ and $\widehat{m}_j^{LL,1,S}$ $(s = A, B)$ are defined according to (6.7)–(6.9). Then there exist random variables $R_{n,i,j}^{LL}(u_j, h, X)$ such that*

$$\widehat{m}_j^{LL,A}(u_j) = \widetilde{m}_j^{LL,A}(u_j) + n^{-1} \sum_{i=1}^{n} R_{n,i,j}^{LL}(u_j, h, X)\varepsilon^i, \tag{6.10a}$$

$$\sup_{0 \le u_j \le 1} \sup_{C_H n^{-1/5} \le h_1,\ldots,h_d \le C'_H n^{-1/5}} |R_{n,i,j}^{LL}(u_j, h, X)| = O_p(1), \tag{6.10b}$$

$$\begin{aligned} &\sup_{0 \le u_j \le 1} \sup_{C_H n^{-1/5} \le h_1,h'_1,\ldots,h_d,h'_d \le C'_H n^{-1/5}} |R_{n,i,j}^{LL}(u_j, h', X) - R_{n,i,j}^{LL}(u_j, h, X)| \\ &\quad = \sum_{j=1}^{d} |h'_j - h_j| O_p(n^{\alpha}) \qquad \textit{for some } \alpha > 0. \end{aligned} \tag{6.10c}$$

*Furthermore, uniformly for $C_H n^{-1/5} \le h_1,\ldots,h_d \le C'_H n^{-1/5}$ and $0 \le u_j \le 1$,*

$$\widehat{m}_j^{LL,A}(u_j) = \widetilde{m}_j^{LL,A}(u_j) + n^{-1} \sum_{i=1}^{n} r_{ij}^{LL}(u_j)\varepsilon^i + o_p(n^{-1/2}), \tag{6.11}$$

$$\widehat{m}_j^{LL,B}(u_j) = m_j(u_j) + O_p(n^{-2/5}), \tag{6.12}$$

*where $r_{ij}^{LL}$ are absolutely uniformly bounded functions that satisfy the Lipschitz condition (6.5). Furthermore, uniformly for $C_H n^{-1/5} \le h_j \le C'_H n^{-1/5}$ and $h_j \le u_j \le 1 - h_j$, we have*

$$\widehat{m}_j^{LL,B}(u_j) = m_j(u_j) + \tfrac{1}{2} m''_j(u_j) h_j^2 \int t^2 K(t)\,dt + o_p(n^{-2/5}). \tag{6.13}$$

6.3. *Proofs.*

PROOF OF THEOREM 6.1. For an additive function $f(x) = f_1(x_1) + \cdots + f_d(x_d)$ we define

$$\widehat{\Psi}_j f(x) = f_1(x_1) + \cdots + f_{j-1}(x_{j-1}) + f_j^*(x_j) + f_{j+1}(x_{j+1}) + \cdots + f_d(x_d),$$

where

$$f_j^*(x_j) = -\sum_{k \neq j} \int f_k(x_k) \frac{\widehat{p}_{jk}(x_j, x_k)}{\widehat{p}_j(x_j)}\, dx_k + \sum_k \int f_k(x_k) \widehat{p}_k(x_k)\, dx_k.$$

According to Lemma 3 in [9], we have for $\widehat{m}^{NW,A}(x) = \widehat{m}_0^{NW,A} + \widehat{m}_1^{NW,A}(x_1) + \cdots + \widehat{m}_d^{NW,A}(x_d)$,

$$\widehat{m}^{NW,A}(x) = \sum_{s=0}^{\infty} \widehat{T}^s \widehat{\tau}(x).$$

Here, $\widehat{T} = \widehat{\Psi}_d \cdots \widehat{\Psi}_1$ and

$$\begin{aligned}\widehat{\tau}(x) = {} & \widehat{\Psi}_d \cdots \widehat{\Psi}_2 [\widetilde{m}_1^{NW,A}(x) - \widetilde{m}_{0,1}^{NW,A}] + \cdots + \widehat{\Psi}_d [\widetilde{m}_{d-1}^{NW,A}(x) - \widetilde{m}_{0,d-1}^{NW,A}] \\ & + \widetilde{m}_d^{NW,A}(x) - \widetilde{m}_{0,d}^{NW,A},\end{aligned}$$

where, in a slight abuse of notation, $\widetilde{m}_j(x) = \widetilde{m}_j(x_j)$ and $\widetilde{m}_{0,j} = \int \widetilde{m}_j(x_j) \times \widehat{p}_j(x_j)\, dx_j$.

We now decompose

$$\widehat{m}^{NW,A}(x) = \widetilde{m}^{NW,A}(x) + \sum_{s=0}^{\infty} \widehat{T}^s (\widehat{\tau} - \widetilde{m}^{NW,A})(x) + \sum_{s=1}^{\infty} \widehat{T}^s \widetilde{m}^{NW,A}(x), \tag{6.14}$$

where $\widetilde{m}^{NW,A}(x) = \widetilde{m}_1^{NW,A}(x_1) + \cdots + \widetilde{m}_d^{NW,A}(x_d)$. We will show that there exist absolutely bounded functions $a^i(x)$ with $|a^i(x) - a^i(y)| \le C\|x - y\|$ for a constant $C$ such that

$$\sum_{s=1}^{\infty} \widehat{T}^s \widetilde{m}^{NW,A}(x) = n^{-1} \sum_{i=1}^{n} a^i(x) \varepsilon^i + o_p(n^{-1/2}) \tag{6.15}$$

uniformly for $C_H n^{-1/5} \le h_j \le C'_H n^{-1/5}$ and $0 \le x_j \le 1$. A similar claim holds for the second term on the right-hand side of (6.14). This immediately implies (6.3).

For the proof of (6.15) we show that there exist absolutely bounded functions $b^i$ with $|b^i(x) - b^i(y)| \le C\|x - y\|$ for a constant $C$ such that

$$\widehat{T} \widetilde{m}^{NW,A}(x) = n^{-1} \sum_{i=1}^{n} b^i(x) \varepsilon^i + o_p(n^{-1/2}), \tag{6.16}$$

$$\sum_{s=1}^{\infty} \widehat{T}^s \widetilde{m}^{NW,A}(x) = \sum_{s=0}^{\infty} T^s \widehat{T} \widetilde{m}^{NW,A}(x) + o_p(n^{-1/2}). \tag{6.17}$$

Here $T = \Psi_d \cdots \Psi_1$ and

$$\Psi_j f(x) = f_1(x_1) + \cdots + f_{j-1}(x_{j-1}) + f_j^{**}(x_j) + f_{j+1}(x_{j+1}) + \cdots + f_d(x_d)$$

for an additive function $f(x) = f_1(x_1) + \cdots + f_d(x_d)$ with

$$f_j^{**}(x_j) = -\sum_{k \neq j} \int f_k(x_k) \frac{p_{jk}(x_j, x_k)}{p_j(x_j)}\, dx_k + \sum_k \int f_k(x_k) p_k(x_k)\, dx_k.$$

Note that (6.15) follows immediately from (6.16) and (6.17), since

$$\begin{aligned}\sum_{s=1}^{\infty} \widehat{T}^s \widetilde{m}^{NW,A}(x) &= \sum_{s=0}^{\infty} T^s \widehat{T} \widetilde{m}^{NW,A}(x) + o_p(n^{-1/2}) \\ &= n^{-1} \sum_{i=1}^{n} \left[ \sum_{s=0}^{\infty} T^s b^i \right](x) \varepsilon^i + o_p(n^{-1/2}).\end{aligned}$$

We prove (6.16) first. For this purpose, one has to consider terms of the form

$$\begin{aligned}S_{kj}(x_j) &= \int \frac{\widehat{p}_{jk}(x_j, x_k)}{\widehat{p}_j(x_j)} \widetilde{m}_k^{NW,A}(x_k)\, dx_k \\ &= n^{-1} \sum_{i=1}^{n} \varepsilon^i \int \frac{\widehat{p}_{jk}(x_j, x_k)}{\widehat{p}_j(x_j) \widehat{p}_k(x_k)} K_{h_k}(x_k, X_k^i)\, dx_k.\end{aligned}$$

We make use of the following well-known facts:

$$\widehat{p}_{jk}(x_j, x_k) = E\{\widehat{p}_{jk}(x_j, x_k)\} + O_p(n^{-3/10} \sqrt{\log n}\,), \tag{6.18}$$

$$\widehat{p}_j(x_j) = E\{\widehat{p}_j(x_j)\} + O_p(n^{-2/5} \sqrt{\log n}\,), \tag{6.19}$$

$$(\partial/\partial x_j) \widehat{p}_{jk}(x_j, x_k) = E\{(\partial/\partial x_j) \widehat{p}_{jk}(x_j, x_k)\} + O_p(n^{-1/10} \sqrt{\log n}\,), \tag{6.20}$$

$$(\partial/\partial x_j) \widehat{p}_j(x_j) = E\{(\partial/\partial x_j) \widehat{p}_j(x_j)\} + O_p(n^{-1/5} \sqrt{\log n}\,), \tag{6.21}$$

uniformly for $C_H n^{-1/5} \le h_j, h_k \le C_H' n^{-1/5}$ and $0 \le x_j, x_k \le 1, \quad 1 \le j,\ k \le d$.

We now argue that

$$\begin{aligned}&S_{kj}(x_j) - n^{-1} \sum_{i=1}^{n} \frac{p_{jk}(x_j, X_k^i)}{p_j(x_j) p_k(X_k^i)} \varepsilon^i \\ &\quad \equiv n^{-1} \sum_{i=1}^{n} \Delta_{kj}(x_j, h_j, h_k) \varepsilon^i = o_p(n^{-1/2}),\end{aligned} \tag{6.22}$$

uniformly in $x_j, h_j, h_k$. From (6.18)–(6.21) and from the expansions of the expectations on the right-hand sides of these equations we get

$$\Delta_{kj}(x_j, h_j, h_k) = O_p(n^{-1/5}),$$

uniformly in $x_j, h_j, h_k$. Furthermore, we have, because of $E[|\varepsilon^i|^{(5/2)+\delta}|X^i] < C$ for some $\delta > 0$, $C < +\infty$, that for a sequence $c_n \to 0$

$$E[|\varepsilon^i|\mathbf{1}(|\varepsilon^i| > n^{2/5})|X^i] \leq c_n n^{-3/5},$$
$$P(|\varepsilon^i| \leq n^{2/5} \text{ for } 1 \leq i \leq n) \to 1.$$

This shows that

$$n^{-1}\sum_{i=1}^n \Delta_{kj}(x_j, h_j, h_k)\varepsilon^i - n^{-1}\sum_{i=1}^n \Delta_{kj}(x_j, h_j, h_k)\varepsilon_*^i = o_p(n^{-1/2}) \tag{6.23}$$

uniformly in $x_j, h_j, h_k$, where

$$\varepsilon_*^i = \varepsilon^i\mathbf{1}(|\varepsilon^i| \leq n^{2/5}) - E[\varepsilon^i\mathbf{1}(|\varepsilon^i| \leq n^{2/5})|X^i].$$

Note now that, with $X = (X^1, \ldots, X^n)$ and $\Delta = n^{1/5}\sup_{k,j,x_j,h_j,h_k}|\Delta_{kj}(x_j, h_j, h_k)|$,

$$\begin{aligned}
P&\left\{n^{-1}\sum_{i=1}^n \Delta_{kj}(x_j, h_j, h_k)\varepsilon_*^i \geq n^{-3/5}|X\right\}\\
&\leq E\left[\exp\left\{n^{-3/10}\sum_{i=1}^n \Delta_{kj}(x_j, h_j, h_k)\varepsilon_*^i\right\}\middle|X\right]\exp(-n^{1/10})\\
&\leq \prod_{i=1}^n E\{\exp\{n^{-3/10}\Delta_{kj}(x_j, h_j, h_k)\varepsilon_*^i\}|X\}\exp(-n^{1/10})\\
&\leq \prod_{i=1}^n [1 + n^{-3/5}\Delta_{kj}^2(x_j, h_j, h_k)E[(\varepsilon_*^i)^2|X^i]\exp\{n^{-3/10}\Delta n^{-1/5}2n^{2/5}\}]\\
&\qquad \times \exp(-n^{1/10})\\
&\leq \exp\left\{\Delta^2 \sup_{1\leq i\leq n} E[(\varepsilon_*^i)^2|X^i]\exp(2n^{-1/10}\Delta)\right\}\exp(-n^{1/10})\\
&\leq M_n\exp(-n^{1/10})
\end{aligned}$$

with a random variable $M_n = O_p(1)$. Together with (6.23) this inequality shows that (6.22) uniformly holds on any grid of values of $x_j, h_j$ and $h_k$ with cardinality being of a polynomial order of $n$. For $\alpha_1, \alpha_2 > 0$ large enough and for a random variable $R_n = O_p(1)$, one can show

$$\begin{aligned}
&|\Delta_{kj}(x_j', h_j', h_k') - \Delta_{kj}(x_j, h_j, h_k)|\\
&\qquad \leq R_n(n^{\alpha_1}|x_j' - x_j| + n^{\alpha_2}|h_j' - h_j| + n^{\alpha_3}|h_k' - h_k|).
\end{aligned}$$

This implies that (6.22) holds uniformly for $0 \leq x_j \leq 1$ and $C_H n^{-1/5} \leq h_j, h_k \leq C_H' n^{-1/5}$. By consideration of other terms similar to $S_{kj}(x_j)$, one may complete the proof of (6.16).

We now come to the proof of (6.17). With the help of (6.18)–(6.21) one can show by using the Cauchy–Schwarz inequality that

$$\sup_{\|f\|\le 1}\sup_{0\le x_1,\dots,x_d\le 1}|\widehat{T}f(x)-Tf(x)|=O_p(n^{-1/10}\sqrt{\log n}\,). \tag{6.24}$$

Here the first supremum runs over all additive functions $f$ with $\int f^2(x)p(x)\,dx \le 1$. (The slow rate is caused by the fact that $\widehat{p}_j$ is inconsistent at $x_j$ in neighborhoods of 0 and 1.) Furthermore, in [9] it has been shown that

$$\sup_{\|f\|\le 1}\sup_{0\le x_1,\dots,x_d\le 1}|\widehat{T}f(x)|=O_p(1), \tag{6.25}$$

$$\sup_{\|f\|\le 1}\|Tf\|<1, \tag{6.26}$$

where $\|Tf\|^2=\int\{Tf(x)\}^2p(x)\,dx$. Claim (6.17) now follows from (6.16), (6.24)–(6.26) and the fact

$$\sum_{s=1}^{\infty}(\widehat{T}^s-T^s)=\sum_{s=1}^{\infty}\sum_{t=0}^{s-1}\widehat{T}^t(\widehat{T}-T)T^{s-1-t}.$$

Proof of (6.2a)–(6.2c). Formula (6.2a) is given by the definition of $\widehat{m}_j^{NW,A}$. Claim (6.2b) follows as in the proof of (6.3). For the proof of (6.2c) one uses bounds on the operator norm of $\widehat{T}_{h'}-\widehat{T}_h$, where $\widehat{T}_h$ is defined as $\widehat{T}$ with bandwidth tuple $h$. □

Proof of (6.4) and (6.6). These claims follow by a slight modification of the arguments used in the proof of Theorem 4 in [9]. There it has been shown that (6.6) holds for fixed bandwidths $h_1,\dots,h_d$ and uniformly for $u_j$ in a closed subinterval of $(0,1)$. The arguments can be easily modified to get uniform convergence for $h_j\le u_j\le 1-h_j$ and $C_Hn^{-1/5}\le h_1,\dots,h_d\le C'_Hn^{-1/5}$. In Theorem 4 in [9] a wrong value was given for $\gamma_{n,j}$; see the wrong proof of (114) in [9]. A correct calculation gives $\gamma_{n,j}$ as stated here. □

The proof of Theorem 6.1 is complete. □

Proof of Theorem 6.2. Theorem 6.2 follows with similar arguments as in the proof of Theorem 6.1. Now one can use Theorem 4′ of [9]. For the proof of (6.13) note that we use another norming for $\widehat{m}_j$ (cf. (6.9) with (52) in [9]). Formula (6.13) follows from Theorem 4′ of [9] by noting that $\int\widehat{m}_j^{LL,1,B}(u_j)\times\widehat{p}_j^1(u_j)\,du_j=-\gamma_{n,j}+o_P(n^{2/5})$ with $\gamma_{n,j}$ defined as in Theorem 4′ of [9]. □

PROOF OF THEOREM 2.1. With $w_i = w(X^i)\mathbf{1}(C_H' n^{-1/5} \le X_j^i \le 1 - C_H' n^{-1/5}$ for $1 \le j \le d)$, we get

$$RSS(h) - ASE(h) = \frac{1}{n}\sum_{i=1}^n w_i(\varepsilon^i)^2 - \frac{2}{n}\sum_{i=1}^n w_i\{\widehat{m}^{NW}(X^i) - m(X^i)\}\varepsilon^i,$$

where $\widehat{m}^{NW}(x) = \widehat{m}_0^{NW} + \widehat{m}_1^{NW}(x_1) + \cdots + \widehat{m}_d^{NW}(x_d)$ and $m(x) = m_0 + m_1(x_1) + \cdots + m_d(x_d)$. We will show that uniformly for $C_H n^{-1/5} \le h_1, \ldots, h_d \le C_H' n^{-1/5}$,

$$\frac{1}{n}\sum_{i=1}^n w_i\{\widehat{m}^{NW,B}(X^i) - m(X^i)\}\varepsilon^i = o_p(n^{-4/5}) \tag{6.27}$$

and

$$\begin{aligned}\frac{1}{n}\sum_{i=1}^n &w_i \widehat{m}^{NW,A}(X^i)\varepsilon^i \\ &= \frac{1}{n}\sum_{i=1}^n w_i(\varepsilon^i)^2 K(0)\sum_{j=1}^d \frac{1}{nh_j} + o_p(n^{-4/5}),\end{aligned} \tag{6.28}$$

where for $S = A, B$ we write

$$\widehat{m}^{NW,S}(x) = \widehat{m}_0^{NW,S} + \widehat{m}_1^{NW,S}(x_1) + \cdots + \widehat{m}_d^{NW,S}(x_d).$$

The statement of Theorem 2.1 immediately follows from (6.27) and (6.28).

For the proof of (6.27) one can proceed similarly as in the proof of (6.22). Note that

$$\sup_{1\le i\le n}\ \sup_{C_H n^{-1/5}\le h_1,\ldots,h_d\le C_H' n^{-1/5}} n^{2/5}|w_i\{\widehat{m}^{NW,B}(X^i) - m(X^i)\}| = O_p(1),$$

and that differences between values of $w^i\{\widehat{m}^{NW,B}(X^i) - m(X^i)\}$ evaluated for different bandwidth tuples $(h_1', \ldots, h_d')$ and $(h_1, \ldots, h_d)$ can be bounded by $\sum_j |h_j' - h_j| O_p(n^\alpha)$ with $\alpha$ large enough.

For the proof of (6.28) we note first that by application of Theorem 6.1,

$$\begin{aligned}\frac{1}{n}\sum_{i=1}^n &w_i\widehat{m}^{NW,A}(X^i)\varepsilon^i \\ &= \frac{1}{n}\sum_{i=1}^n w_i\widetilde{m}^{NW,A}(X^i)\varepsilon^i + \frac{1}{n^2}\sum_{i,k=1}^n\sum_{j=1}^d w_i R_{n,k,j}(X_j^i, h, X)\varepsilon^i\varepsilon^k \\ &= \frac{1}{n^2}\sum_{i=1}^n\sum_{j=1}^d w_i h_j^{-1} K(0)(\varepsilon^i)^2\end{aligned}$$

$$+\frac{1}{n^2}\sum_{i\neq k}\sum_{j=1}^{d} w_i h_j^{-1} K[h_j^{-1}(X_j^i - X_j^k)]\varepsilon^i\varepsilon^k$$

$$+\frac{1}{n^2}\sum_{i=1}^{n}\sum_{j=1}^{d} w_i R_{n,i,j}(X_j^i, h, X)(\varepsilon^i)^2$$

$$+\frac{1}{n^2}\sum_{i\neq k}\sum_{j=1}^{d} w_i R_{n,k,j}(X_j^i, h, X)\varepsilon^i\varepsilon^k$$

$$= T_1(h) + \cdots + T_4(h).$$

Now, it is easy to check that uniformly for $C_H n^{-1/5} \le h_1, \ldots, h_d \le C_H' n^{-1/5}$,

$$T_1(h) = \frac{1}{n}\sum_{i=1}^{n} w_i(\varepsilon^i)^2 K(0)\sum_{j=1}^{d}\frac{1}{nh_j}\{1 + O_p(n^{-1/2}\sqrt{\log n}\,)\},$$

$$|T_3(h)| \le O_p(1)\frac{1}{n^2}\sum_{i=1}^{n}(\varepsilon^i)^2 = O_p(n^{-1}).$$

So, it remains to show

$$T_2(h) = o_p(n^{-4/5}), \tag{6.29}$$

$$T_4(h) = o_p(n^{-4/5}). \tag{6.30}$$

We will show (6.29). Claim (6.30) follows by slightly simpler arguments. For (6.29) it suffices to show that for $1 \le j \le d$

$$T_{2,j}^*(h) \equiv \frac{1}{n^2}\sum_{i\neq k} w_i h_j^{-1} K[h_j^{-1}(X_j^i - X_i^k)]\eta^i\eta^k = o_p(n^{-4/5}), \tag{6.31}$$

where $\eta^i = \varepsilon^i \mathbf{1}(|\varepsilon^i| \le n^\alpha) - E[\varepsilon^i \mathbf{1}(|\varepsilon^i| \le n^\alpha)|X^i]$ with $1/\gamma < \alpha < 1/4$. The constant $\gamma$ was introduced in assumption (A3). It holds that $E|\varepsilon^i|^\gamma < C'$ for some $C' < \infty$; see assumption (A3). Note that

$$P(|\varepsilon^i| > n^\alpha \text{ for some } i \text{ with } 1 \le i \le n) \le nE|\varepsilon^1|^\gamma n^{-\alpha\gamma} \to 0,$$

$$E[|\varepsilon^i|\mathbf{1}(|\varepsilon^i| > n^\alpha)|X^i] \le E[|\varepsilon^i|^\gamma|X^i]n^{-(\gamma-1)\alpha} \le C' n^{-(\gamma-1)\alpha} = O(n^{-3/4}).$$

We apply an exponential inequality for $U$-statistics. Let

$$\kappa_n^2 = E\{2^{-1}(w_i + w_k)n^{1/10}K[h_j^{-1}(X_j^i - X_j^k)]\eta^i\eta^k\}^2,$$

$$M_n = \sup\{2^{-1}(w_i + w_k)n^{1/10}K[h_j^{-1}(X_j^i - X_j^k)]\eta^i\eta^k\},$$

where the supremum in the definition of $M_n$ is over the whole probability space. We note that $\kappa_n^2 = O(1)$ and $M_n$ is bounded by a constant which is $O(n^{2\alpha}n^{1/10})$. According to Theorem 4.1.12 in [2], for constants $c_1, c_2 > 0$ and $0 < \delta < \frac{1}{2} - 2\alpha$,

$$\begin{aligned}
&P(|T^*_{2,j}(h)| \geq n^{-4/5-\delta}) \\
&\qquad \leq P\left(\left|n^{-1}\sum_{i\neq k} w_i n^{1/10} K[h_j^{-1}(X_j^i - X_j^k)]\eta^i\eta^k\right| \geq C_H n^{1/10-\delta}\right) \\
&\qquad \leq c_1 \exp\left[-\frac{c_2 n^{1/10-\delta}}{\kappa_n + \{M_n n^{(1/10-\delta)/2} n^{-1/2}\}^{2/3}}\right].
\end{aligned}$$

This gives with $\rho = (1 - 2\delta - 4\alpha)/3 > 0$ and a constant $c_3 > 0$,

$$P(|T^*_{2,j}(h)| \geq n^{-4/5-\delta}) \leq c_1 \exp(-c_3 n^\rho).$$

Together with $|T^*_{2,j}(h') - T^*_{2,j}(h)| \leq cn^\alpha |h'_j - h_j|$ for $c, \alpha > 0$ large enough, this implies (6.31). □

PROOF OF THEOREM 2.2. Claim (2.11a) follows from the expansions of Theorem 6.1. For the proof of (2.11b) note that

$$RSS(h) = n^{-1}\sum_{i=1}^n w(X^i)(\varepsilon^i)^2 + o_p(1)$$

because of (2.11a) and Theorem 2.1. □

PROOF OF THEOREMS 3.1–3.3. Theorems 3.1 and 3.2 follow with similar arguments as in the proofs of Theorems 2.1 and 2.2. For the proof of Theorem 3.3, one uses

$$\sup_{0\leq v_j\leq 1} |\widehat{m}_j^{LL}(v_j) - m_j(v_j)| = O_p(n^{-2/5}\sqrt{\log n}\,). \tag{6.32}$$

This can be shown by use of the expansions of Theorem 6.2. By standard arguments in local polynomial regression (see [3], e.g.), it follows that for $\widehat{m}''_j(u_j)$ defined in (3.16),

$$\begin{aligned}
&\widehat{m}''_j(u_j) - m''_j(u_j) \\
&\qquad = 2g_j^{-3}\int L^*[g_j^{-1}(v_j - u_j)] \\
&\qquad\qquad\qquad \times \{\widehat{m}_j^{LL}(v_j) - m_j(u_j) - m'_j(u_j)(v_j - u_j) \\
&\qquad\qquad\qquad\qquad - \tfrac{1}{2}m''_j(u_j)(v_j - u_j)^2\}\,dv_j,
\end{aligned}$$

where $L^*$ is the so-called *equivalent kernel* having the property that $\int L^*(v_j)\,dv_j = \int v_j L^*(v_j)\,dv_j = 0$ and $\int v_j^2 L^*(v_j)\,dv_j = 1$. Application of (6.32) gives

$$g_j^{-3}\int L^*[g_j^{-1}(v_j - u_j)]\{\widehat{m}_j^{LL}(v_j) - m_j(v_j)\}\,dv_j = o_p(1).$$

Now, the fact that the function itself and its first two derivatives at $u_j$ of $\nu(\cdot) \equiv m_j(\cdot) - m_j(u_j) - m_j'(u_j)(\cdot - u_j) - m_j''(u_j)(\cdot - u_j)^2/2$ are all zero yields

$$g_j^{-3} \int L^*[g_j^{-1}(v_j - u_j)] \times \{m_j(v_j) - m_j(u_j) - m_j'(u_j)(v_j - u_j) - \tfrac{1}{2} m_j''(u_j)(v_j - u_j)^2\}\, dv_j = o(1). \qquad \square$$

**Acknowledgments.** We are grateful for the helpful and constructive comments of four reviewers. Assistance of Jun Hyeok Hwang for the numerical work is greatly appreciated.

Department of Economics
University of Mannheim
L 7, 3-5
68131 Mannheim
Germany
e-mail: emammen@rumms.uni-mannheim.de

Department of Statistics
Seoul National University
Seoul 151-747
Korea
e-mail: bupark@stats.snu.ac.kr